\documentclass[12pt,reqno]{article}%
\usepackage{latexsym}
\usepackage{amssymb}
\usepackage{amsmath}
\usepackage[utf8]{inputenc}
\usepackage{amsfonts}
\usepackage{graphicx}%
\providecommand{\U}[1]{\protect\rule{.1in}{.1in}}
\numberwithin{equation}{section}
\newtheorem{theorem}{Theorem}[section]

\newtheorem{definition}[theorem]{Definition}

\newtheorem{lemma}[theorem]{Lemma}

\newtheorem{remark}[theorem]{Remark}

\begin{document}

\title{Boundary Behavior of Subelliptic Parabolic Equations on Time-Dependent Domains}
\author{Elin G{\"o}tmark\thanks{email: elin.gotmark@math.umu.se}\linebreak\\Department of Mathematics\linebreak\\Ume{\aa } University\linebreak\\S-90187 Ume{\aa } \\Sweden
\and Marie Frentz\thanks{email: e.m.frentz@lse.ac.uk}\\Department of Mathematics\\London School of Economics\\Houghton Street\\London\\WC2A 2AE\\United Kingdom}
\maketitle

\begin{abstract}
\noindent In this paper we study the boundary behavior of solutions of a
divergence-form subelliptic heat equation in a time-varying domain
$\Omega\subset\mathbb{R}^{n+1}$, structured on a set of vector fields
$X=\{X_{1},...,X_{m}\}$ with $C^{\infty}$-coefficients satisfying
H{\"{o}}rmander's finite rank condition. Assuming that $\Omega$ is an $X$-NTA
domain, we first prove a Dahlberg type estimate comparing $\omega=\omega
_{X}^{(x,t)}$, which is the $X$-caloric measure at $(x,t)\in\Omega_{T}$, and
the Green function of $H$. We then prove a backward Harnack inequality, the
doubling property for $\omega_{X}$, the H{\"{o}}lder continuity at the
boundary for quotients of solutions of $H$, and a Fatou theorem. \newline

\noindent

\noindent\emph{2010 Mathematics Subject Classification:} 35K70, 31B25.
\newline\noindent

\noindent\textit{Keywords and phrases: Boundary behavior, subelliptic
parabolic, divergence form, caloric measure, time-dependent domains.}

\end{abstract}

\section{Introduction}

In this paper, we study the boundary behavior of solutions of the subelliptic
heat operator
\begin{equation}
H=L-\partial_{t},\mbox{ where }L=\sum_{i,j=1}^{m}X_{i}^{\ast}(a_{ij}%
(x,t)X_{j}u),\ (x,t)\in\mathbb{R}^{n}\times\mathbb{R}. \label{1.1}%
\end{equation}
The system $X=\{X_{1},...,X_{m}\},$ $m<n$, is a set of vector fields with
$C^{\infty}$-coefficients satisfying H{\"{o}}rmander's finite rank condition,
i.e.,
\begin{equation}
\mbox{ the rank of $\mbox{Lie}[X_1,...,X_m]$ equals $n$}. \label{1.1x}%
\end{equation}
Moreover, $X_{i}^{\ast}$ denotes the adjoint operator of $X_{i}$, which for
$X_{i}=\sum b_{j}^{i}(x)\partial_{j}$ is given by $X_{i}^{\ast}=-X_{i}%
-\sum\frac{\partial b_{j}^{i}(x)}{\partial x_{j}}$. We let $\hat{H}$ denote
the adjoint subelliptic heat operator.

Let $\Omega\subset\mathbb{R}_{x}^{n}\times\mathbb{R}_{t}$ be a domain (that
is, an open and connected set) which is bounded in the $\mathbb{R}_{x}^{n}$
variables, where the subscript $x$ (respectively $t$) indicates that we are
only concerned with the space (respectively time) variables. We will further
assume that $\Omega$ is contained in $\{(x,t):t>0\}$ and that $(0,0)$ lies on
the boundary of $\Omega$. The metric on $\Omega$ will be $d_{p}%
(x,t,y,s)=(d(x,y)^{2}+|t-s|)^{1/2}$, where $d(x,y)$ is the
Carnot-Carath\'{e}odory metric on $\mathbb{R}^{n}$ induced by the vector
fields $X$; see (\ref{storm}) for a definition. We will assume that $\Omega$
is an $X$-NTA domain with constants $r_{0}$ and $M$, see Definition
\ref{defNTA}. In the first sections, Section \ref{wrangle}-Section \ref{SecH},
we will work on the bounded domain $\Omega_{T}=\Omega\cap\{0\leq t<T\}$.

We assume that the matrix $A=(a_{ij}(x,t))$ in the operator (\ref{1.1}) is
symmetric, bounded, and uniformly elliptic, i.e., that there exists
$\lambda>0$ such that
\[
\lambda|\xi|^{2}\leq\sum_{i,j}^{m}a_{ij}(x,t)\xi_{i}\xi_{j}\leq\lambda
^{-1}|\xi|^{2}.
\]
Note that uniform ellipticity is only required in $m$ out of $n$ spatial
directions, $m<n$, in contrast to the classical case where $m=n$ and
$X_{i}=\partial x_{i}$. We also assume that the coefficients $a_{ij}$ are
smooth. A natural relaxation would have been to assume that the coefficients
$a_{ij}$ are merely $\alpha$-H\"{o}lder continuous functions with respect to
the metric $d_{p}$. However, at present, we are unable to derive such results
and at the end of this section we shortly explain why.

Let $\partial_{p}\Omega_{T}$ be the parabolic boundary of $\Omega_{T}$, that
is, $\{(x,t)\in\partial\Omega_{T}:t<T\}$. By using the maximum principle of
Bony \cite{B} one can show that there exists, for any $f$ which is continuous
on $\partial_{p}\Omega_{T}$, a unique Perron-Wiener-Brelot solution
$u_{f}^{\Omega_{T}}$ to the Dirichlet problem
\[
Hu=0\mbox{ in $\Omega_T$, $u=f$ on $\partial_p\Omega_T$}.
\]
Using that $\Omega$ is an $X$-NTA domain, one can also prove that
$\partial_{p}\Omega_{T}$ is regular for the Dirichlet problem for $H$, that
is, $u_{f}^{\Omega_{T}}$ is continuous on the closure of $\Omega_{T}$ and
achieves the boundary data continuously.

The maximum principle and the Riesz representation theorem yield the
representation formula
\begin{equation}
u(x,t)=\int\limits_{\partial_{p}\Omega_{T}}f(y,s)d\omega
(x,t,y,s)\mbox{ whenever
}(x,t)\in\Omega_{T} \label{caloricM}%
\end{equation}
where $\omega(x,t)=\omega_{X}(x,t)$ is the $X$-caloric or $X$-parabolic
measure at $(x,t)$ associated to the subelliptic heat operator $H$. If
$E\subset\partial_{p}\Omega_{T}$ is a Borel set, then $\omega(x,t,E)$, or the
measure of $E$ with respect to $\omega(x,t)$, is the solution of
\[
Hu=0\mbox{ in $\Omega_T$, $u=\chi_E$ on $\partial_p\Omega_T$},
\]
where $\chi_{E}$ is the characteristic function of the set $E$. Sometimes we
will wish to construct the measure $\omega(x,t)$ with respect to another set
than $\Omega_{T}$, say, $F\subset\mathbb{R}_{x}^{n}\times\mathbb{R}_{t}$. In
that case, we will denote the measure by $\omega(x,t,F)$, and the measure of
$E\subset\partial_{p}F$ will be denoted $\omega(x,t,E,F)$. By $\hat{\omega}$
we mean the $X$-parabolic measure associated to the adjoint operator $\hat{H}$.

We let $S_{T}=\{(x,t)\in\partial\Omega:t\in(0,T)\}$ and define the following
sets:
\begin{align*}
B_{d}(x,r)  &  =\{y\in\mathbb{R}^{n}:\ d(x,y)<r\}\\
C_{r}(x,t)  &  =B_{d}(x,r)\times(t-r^{2},t+r^{2})\\
\Delta(x,t,r)  &  =C_{r}(x,t)\cap\partial_{p}\Omega_{T},
\end{align*}
where the last set is defined whenever $(x,t)\in\partial_{p}\Omega_{T}$. When
we say that a constant $c$ depends on $H$, we mean that it depends on $n $,
$m$, the ellipticity constant $\lambda$, the vector fields $X$, and the
$\alpha$-H\"{o}lder norms of the $a_{ij}$.

Let $G(x,t,y,s)$ be the Green function with respect to the operator $H$ in
(\ref{1.1}) and the set $\Omega_{T}$. Our first result is a Dahlberg type
estimate comparing the $X$-caloric measure and the Green function for
$\Omega_{T}$ and $H$:

\begin{theorem}
\label{lem4.10}\emph{(Dahlberg estimate)} Let $X=\{X_{1},...,X_{m}\}$ be a
system of smooth vector fields satisfying \eqref{1.1x}. Let $\Omega
\subset\mathbb{R}^{n+1}$ be a $X$-NTA-domain in the sense of Definition
\ref{defNTA} with parameters $M$, $r_{0}$. Let $(x_{0},t_{0})\in S_{T}$ and
$r<\min\{r_{0}/2,\sqrt{(T-t_{0})/4},\sqrt{t_{0}/4}\}$. Assume that
$(x,t)\in\Omega_{T}$ and that $|t-t_{0}|\geq5r^{2}$.Then there exists
$c=c(H,M,r_{0})$, $1\leq c<\infty$, such that
\begin{align*}
c^{-1}\ |B_{d}(x,r)|\ G(x,t,A_{r}^{+}(x_{0},t_{0}))  &  \leq\omega
(x,t,\Delta(x_{0},t_{0},r/2))\leq\\
&  \leq c\ |B_{d}(x,r)|\ G(x,t,A_{r}^{-}(x_{0},t_{0}))\mbox{ if }t\geq t_{0};
\end{align*}%
\begin{align*}
c^{-1\ }|B_{d}(x,r)|\ G(A_{r}^{-}(x_{0},t_{0}),x,t)  &  \leq\hat{\omega
}(x,t,\Delta(x_{0},t_{0},r/2))\leq\\
&  \leq c\ |B_{d}(x,r)|\ G(A_{r}^{+}(x_{0},t_{0}),x,t)\mbox{ if }t\leq t_{0}.
\end{align*}

\end{theorem}

\noindent Theorem \ref{lem4.10} is used to prove the doubling condition for
the $X$-caloric measure:

\begin{theorem}
\label{lem4.13}\emph{(Doubling property)} Let $X=\{X_{1},...,X_{m}\}$ be a
system of smooth vector fields satisfying \eqref{1.1x}. Let $\Omega
\subset\mathbb{R}^{n+1}$ be a $X$-NTA-domain in the sense of Definition
\ref{defNTA} with parameters $M$, $r_{0}$. Let $0<\delta\ll1$ be a fixed
constant and let $(x_{0},t_{0})\in S_{T}$ with $\delta^{2}\leq t_{0}\leq
T-\delta^{2}$. Then there exists a constant $c=c(H,M,r_{0},\mbox{diam}(\Omega
),T,\delta)$, $1\leq c<\infty$, such that if $r<\min\{r_{0}/2,\sqrt
{(T-t_{0}-\delta^{2})/4},\sqrt{(t_{0}-\delta^{2})/4}\}$, then
\[
\omega^{\ast}(x,t,\Delta(x_{0},t_{0},r))\leq c\omega^{\ast}(x,t,\Delta
(x_{0},t_{0},r/2)),
\]
where $\omega^{\ast}=\omega$ when $t-t_{0}\geq10r^{2}$ while $\omega^{\ast
}=\hat{\omega}$ when $t_{0}-t\geq10r^{2}$.
\end{theorem}

\noindent Next, we prove the H{\"o}lder continuity of the quotient of two solutions:

\begin{theorem}
\label{holder}\emph{(H\"{o}lder continuity of quotients of solutions)} Let
$X=\{X_{1},...,X_{m}\}$ be a system of smooth vector fields satisfying
\eqref{1.1x}. Let $\Omega\subset\mathbb{R}^{n+1}$ be a $X$-NTA-domain in the
sense of Definition \ref{defNTA} with parameters $M$, $r_{0}$. Let $u$, $v$ be
non-negative solutions to either the subelliptic heat equation, $Hu=0$, or the
adjoint subelliptic heat equation, $\hat{H}u=0$, in $\Omega_{T}\cap
C_{4r}(x_{0},t_{0})$, where $r<\min\{r_{0}/2,\sqrt{(T-t_{0}-\delta^{2}%
)/4},\sqrt{(t_{0}-\delta^{2})/4}\}$. Assume that $u$, $v$ vanish continuously
on $\Delta(x_{0},t_{0},2r)$. Then $u/v$ is H{\"{o}}lder continuous on the
closure of $\Omega_{T}\cap C_{r}(x_{0},t_{0})$.
\end{theorem}

\noindent Lastly, we prove a Fatou theorem. We define a non-tangential region
at $P \in S_{T}$ as
\[
\Gamma_{\alpha}(P) = \{ (x,t) \in\Omega: d_{p}(x,t,P) \leq(1+ \alpha)
d_{p}(x,t,\partial_{p} \Omega) \},
\]
and the non-tangential maximal function as $N_{\alpha}(u)(P) = \sup_{(x,t)
\in\Gamma_{\alpha}(P)} |u(x,t)|$.

\begin{theorem}
\label{fatou}\emph{(Fatou theorem)} If $u$ is a non-negative solution of
$Hu=0$ in $\Omega_{T}$, then the non-tangential limit $N_{\alpha}(u)$ exists
almost everywhere on $S_{T}$ with respect to the measure $\omega$.
\end{theorem}

Fatou type theorems have a long and rich history. In particular, since the
appearance of Fatou theorems in the papers \cite{HW1}, \cite{HW2}\ of Hunt and
Wheeden on non-tangential convergence of harmonic functions in Lipschitz
domains, Fatou type theorems have been an interesting area of study. Later on,
results for second order elliptic equations was obtained by Caffarelli et al.
in \cite{Caf} on Lipschitz domain and by Jerison and Kenig in \cite{JK} on
NTA-domains. It was not until 1998, some thirty years later, corresponding
results was proved in the subelliptic setting by Capogna and Garofalo in
\cite{CG}. When $m=n$ and $X_{i}=\partial_{x_{i}}$, we get classical parabolic
operators, and we refer to the work of Fabes, Garofalo and Salsa, \cite{FGS},
where similar results for parabolic divergence form operators was proved. The
importance of Fatou type theorems stems, amongst others, from the study of
free boundaries. For instance, the results in \cite{CG} were crucial in the
proof of the regularity of the free boundary for subelliptic obstacle problems
in \cite{Dan}. Actually, to advance the study of parabolic subelliptic
obstacle problems initiated in \cite{F} and \cite{FGN} was our main motivation
for carrying out this study.

In \cite{FGGMN} we proved a backward Harnack inequality, Theorem
\ref{lem4.13}, and Theorem \ref{holder} for parabolic sub-elliptic operators
in non-divergence form on a domain $\Omega_{T}=\Omega\times\lbrack0,T)$.
Results from this paper cannot apply verbatim to our situation because of our
time-dependent domain, but several of the proofs in our paper are similar
enough that we refer to that paper. Another difference is that the proofs in
\cite{FGGMN} do not use the Green function. Finally, in \cite{M}, Munive
proves some results similar to ours, on cylindrical (not time-dependent)
domains $\Omega\times(0,T)$ when $L$ in (\ref{1.1}) is given by $L=\sum
_{i,j=1}^{m}X_{i}^{\ast}X_{j}$. Note, in particular, that no attempts are made
to prove Fatou type theorems in \cite{M}.

We remark that, to our surprise, we were not able to assume that the
coefficients $a_{ij}$ are merely H\"{o}lder continuous, but rather we had to
assume that the $a_{ij}$'s are smooth. The problems we encountered, which
finally forced us to make this restriction, were \ (a) the lack of a strong
maximum principle, (b) estimates on fundamental solutions and \ (c) the
question of which domains that are regular for the Dirichlet problem. However,
should these results be available, our results carries over directly to this
more general setting.

The paper is organized as follows. In Section 2 we provide the reader with
necessary background information. In Section \ref{wrangle} we prove a
backward-in-time Harnack inequality at the boundary, see Theorem
\ref{lem4.12}. In Section \ref{SecDahl} we continue our study and prove the
Dahlberg estimate in Theorem \ref{lem4.10}. In Section \ref{SecH} our main
result is Theorem \ref{lem4.11}, from which the H\"{o}lder continuity of
quotients, Theorem \ref{holder}, follows. Finally, in Section \ref{SecFatou}
we prove a Fatou type theorem, Theorem \ref{fatou}.

\setcounter{equation}{0} \setcounter{theorem}{0}

\section{Preliminaries}

\subsection{Definitions and notation \label{prelim}}

We now define the Carnot-Carath\'{e}odory distance between $x,y\in
\mathbb{R}^{n}$, induced by $\{X_{1},...,X_{m}\}$. First, let a piecewise
continuous curve $\gamma:[0,l]\rightarrow\mathbb{R}^{n}$ belong to
$\mathcal{S}(x,y)$ if $\gamma(0)=x$ and $\gamma(l)=y$ and
\[
\langle\gamma^{\prime}(t),\xi\rangle^{2}\leq\sum_{1}^{m}\langle X_{j}%
(\gamma(t)),\xi\rangle^{2}%
\]
for every $\xi\in\mathbb{R}^{n}$; note that this implies that $\gamma^{\prime
}(t)\in\text{Span}(X_{j}(\gamma(t)))$. Moreover, we let $l_{s}(\gamma):=l$ be
called the sub-unitary length of $\gamma$. Then the Carnot-Carath\'{e}odory
distance between $x,y\in\mathbb{R}^{n}$ is
\[
d(x,y)=\inf\{l_{s}(\gamma):\gamma\in\mathcal{S}(x,y)\}.
\]
From \cite[Proposition 1.1]{NSW} we have the following fact: there exist
$C,\epsilon>0$ such that
\[
C|x-y|\leq d(x,y)\leq C^{-1}|x-y|^{\epsilon}%
\]
for all $x,y\in\mathbb{R}^{n}$. If $x\in\mathbb{R}^{n}$, then we define
\[
B_{d}(x_{0},r):=\{x\in\mathbb{R}^{n}:d(x,x_{0})<r\}.
\]
Note that for large $r$ the closure of Carnot-Carath\'{e}odory balls
$B_{d}(x,r)$ may fail to be compact, see \cite[p. 1086]{GN}. But by Chow's
accessibility theorem, $(\mathbb{R}^{n},d)$ is at least a locally compact
space, see \cite{C}, and we can find $R_{0}>0$ such that the closure of any
ball $B_{d}$ contained in $B_{d}(0,R_{0})$ is compact. We will always assume
that $\Omega\subset B_{d}(0,R_{0}/2)\times(0,\infty)$. Also, again by
\cite[Section \ 3]{NSW}, there exists a function $\Lambda(x,r)$ which is a
polynomial in $r,$ with coefficients that depend on $x$, where the terms may
have degrees between $n$ and $Q$ (the so-called local homogeneous dimension of
$\Omega$). We also have the inequality
\[
C\Lambda(x,r)\leq|B_{d}(x,r)|\leq C^{-1}\Lambda(x,r)
\]
for some constant $C$. It follows that
\begin{equation}
|B_{d}(x,ar)|\leq a^{Q}C^{-1}|B_{d}(x,r)| \label{axel}%
\end{equation}
for $a>1$ and $B_{d}(x,ar)\subset B_{d}(0,R_{0})$. Since $\Omega$ is compact,
it also follows that there exists a constant $C$ depending on $\Omega$, such
that
\begin{equation}
Cr^{Q}\leq|B_{d}(x,r)|\leq C^{-1}r^{n}. \label{motaggsvamp}%
\end{equation}

Now, let $d_{p}(x,t,y,s)=(d(x,y)^{2}+|t-s|)^{1/2}$; then $d_{p}$ is a metric
on $\Omega$. We define
\[
B_{d}^{t_{0}}(x_{0},r):= \{(x,t) \in\mathbb{R}^{n+1}: d(x,x_{0}) <r,
t=t_{0}\}.
\]
By $|B_{d}^{t_{0}}(x_{0},r)|$ we mean the Lebesgue measure of $B_{d}(x_{0},r)
\subset\mathbb{R}^{n}$. If $(x,t)\in\mathbb{R}^{n+1}$, $(x_{0},t_{0}) \in
S_{T}$ and $r,\rho>0$, we let
\[
C_{r,\rho}(x,t):=B_{d}(x,r)\times(t-\rho^{2},t+\rho^{2})
\]
\[
\Delta(x_{0},t_{0},r,\rho):= C_{r,\rho}(x_{0},t_{0}) \cap\partial_{p}
\Omega_{T}.
\]
A cylinder $C_{r,\rho}(x,t)$ is called $(M,X)$-non-tangential in $\Omega_{T}$
if
\[
M^{-1}<\frac{r}{d(C_{r,\rho}(x,t),\partial\Omega_{T})}<M.
\]

\begin{definition}
\label{def1.5} Assume that $(x,t),(y,s)\in\Omega_{T}$ satisfy $(s-t)^{1/2}%
\geq\eta^{-1}d_{p}(x,t,y,s)$ for some $\eta>1$. A sequence of cylinders in
$\Omega_{T}$, $C_{r_{1},\rho_{1}}(x_{1},t_{1})$,..., $C_{r_{l},\rho_{l}}%
(x_{l},t_{l})$, will be called a parabolic Harnack chain of length $l$ joining
$(x,t)$ to $(y,s)$ if there exists a constant $c=c(\nu)$ such that
\begin{align*}
(i)  &  c(\eta)^{-1}\leq\frac{\rho_{i}}{r_{i}}\leq c(\eta
)\mbox{ for }i=1,2,\dots,l\\
(ii)  &  t_{i+1}-t_{i}\geq c(\eta)^{-1}r_{i}^{2},\mbox{ for }i=1,2,\dots,l\\
(iii)  &  \mbox{$C_{r_i,\rho_i}(x_i,t_i)$ is $(M,X)$-non-tangential
in $\Omega_T$  for  $i = 1, 2, \dots, l$}\\
(iv)  &  (x,t)\in C_{r_{1},\rho_{1}}(x_{1},t_{1}),\ (y,s)\in C_{r_{l},\rho
_{l}}(x_{l},t_{l})\\
(v)  &  C_{r_{i+1},\rho_{i+1}}(x_{i+1},t_{i+1})\cap C_{r_{i},\rho_{i}}%
(x_{i},t_{i})\neq\emptyset\mbox{ for }i=1,2,\dots,l-1.
\end{align*}

\end{definition}

We will assume that $\Omega$ is a $X$-NTA-domain in the following sense:

\begin{definition}
\label{defNTA} Let $X=\{X_{1},...,X_{m}\}$ be a system of smooth vector fields
satisfying \eqref{1.1x}. We say that the domain $\Omega\subset\mathbb{R}%
^{n+1}$ defined in the introduction is a non-tangentially accessible domain
with respect to the system $X=\{X_{1},...,X_{m}\}$, in the following referred
to as an $X$-NTA-domain, if there exists $M\geq1$ and $r_{0}>0$ such that the
following are fulfilled:

\begin{enumerate}
\item[(i)] Given $(x_{0},t_{0})\in S_{T}$ and $r\leq\min(r_{0},\sqrt{t_{0}%
/2})$, there exists $A_{r}^{+}(x_{0},t_{0})=(y,s)\in\Omega$ with
$s=t_{0}+2r^{2}$ and $M^{-1}r<d(x_{0},y)<r$ and $d_{p}(A_{r}^{+}(x_{0},t_{0}),
\partial_{p} \Omega)\geq M^{-1}r$. Similarly, there also exists $A_{r}%
^{-}(x_{0},t_{0})$ with the time coordinate equal to $t_{0} -2r^{2}$, and
$A_{r}(x_{0},t_{0})$ with time coordinate equal to $t_{0}$.

\item[(ii)] Condition (i) with $\Omega$ replaced by $\Omega^{C}$. We call the
corresponding points $\underline{A}_{r}^{+}(x_{0},t_{0})$, $\underline{A}
_{r}^{-}(x_{0},t_{0})$, and $\underline{A}_{r}(x_{0},t_{0})$.

\item[(iii)] If $(x,t),(y,s)\in\Omega$ satisfy $( s - t )^{1/2}\geq\eta^{ - 1}
d_{p}(x,t,y,s)$ for some $\eta>1$, and $d(x,\partial\Omega)>\epsilon$,
$d(y,\partial\Omega)>\epsilon$, $(T-s)>\epsilon^{2}$, $t>\epsilon^{2}$ and
$d_{p}(x, t,y, s)<c\epsilon$ for some $\epsilon>0$, then there exists a
parabolic Harnack chain of length $l$ joining $(x,t)$ to $(y,s)$, where $l$
can be chosen independently of $\epsilon$ but depending on $\eta$ and $c$.
\end{enumerate}
\end{definition}

Finally, we need to define H\"{o}lder continuity with respect to the
Carnot-Carath\'{e}odory metric: Let $U\subset\mathbf{R}^{n+1}$ be a bounded
domain and let $\alpha\in(0,1]$. Given $U$ and $\alpha$ we define the
{H\"{o}lder space $C^{0,\alpha}(U)$ as $C^{0,\alpha}(U)=\{u:U\rightarrow
\mathbf{R}:\ ||u||_{C^{0,\alpha}(U)}<\infty\}$, where
\begin{align}
||u||_{C^{0,\alpha}(U)}  &  =|u|_{C^{0,\alpha}(U)}+||u||_{L^{\infty}%
(U)},\label{storm}\\
|u|_{C^{0,\alpha}(U)}  &  =\sup\left\{  \frac{|u(x,t)-u(y,t)|}{d_{p}%
((x,t),(y,s))^{\alpha}}:(x,t),(y,t)\in U,\ (x,t)\neq(y,s)\right\}  .\nonumber
\end{align}
Given a multiindex $I=(i_{1},i_{2},...,i_{m})$, with $1\leq i_{j}\leq m$, we
define $|I|=m$ and $X^{I}u=X_{i_{1}}X_{i_{2}}\cdots X_{i_{m}}u$. Given $U$,
$\alpha$ and an arbitrary non-negative integer $k$ we let $C^{k,\alpha
}(U)=\{u:U\rightarrow\mathbf{R}:\ ||u||_{C^{k,\alpha}(U)}<\infty\}$, where
\[
||u||_{C^{k,\alpha}(U)}=\sum_{|I|+2h\leq k}||\partial_{t}^{h}X^{I}%
u||_{C^{0,\alpha}(U)}.
\]
We also define the class $\Gamma^{2}(U)$ to be the set of all continuous
functions $u$ on $U$ such that $\partial_{t}u$ as well as $X_{i}u$ and
$X_{i}X_{j}u$ are continuous on $U$ for all $i$ and $j$. }

\subsection{The Dirichlet problem}

We will study the Dirichlet problem
\begin{equation}
Hu=g\mbox{ in $\Omega_T$, $u=f$ on $\partial_p \Omega_T$}, \label{dirp}%
\end{equation}
where $f\in C(\partial_{p}\Omega_{T})$, with $C(\partial_{p}\Omega_{T})$
denoting the space of real-valued functions continuous on $\partial_{p}%
\Omega_{T}$. First, we have a strong maximum principle which follows from
\cite[Theoreme 3.2]{B}:

\begin{theorem}
\label{strongmax}\emph{(Strong maximum principle)} Let $X=\{X_{1},...,X_{m}\}
$ be a system of smooth vector fields satisfying \eqref{1.1x}, and let
$\Omega_{T}\subset\mathbb{R}^{n}$ be a bounded domain. Assume that $u\in
\Gamma^{2}(\Omega_{T})$ and that $u\leq0$ in $\Omega_{T}.$Then the following hold:

\begin{enumerate}
\item If $Hu\geq0$ in $\Omega_{T}$ and if $u(x_{0},t_{0})=0$ for some
$(x_{0},t_{0})\in\Omega_{T}$, then $u(x,t)\equiv0$ whenever $(x,t)\in
\Omega_{T}\cap\{t:t\leq t_{0}\}$.

\item If $\hat{H}u\geq0$ in $\Omega_{T}$ and if $u(x_{0},t_{0})=0$ for some
$(x_{0},t_{0})\in\Omega_{T}$, then $u(x,t)\equiv0$ whenever $(x,t)\in
\Omega_{T}\cap\{t:t\geq t_{0}\}$.
\end{enumerate}
\end{theorem}

The following theorem shows that we can solve the Dirichlet problem on $X$-NTA
domains and that all points in $\partial_{p}\Omega_{T}$ are regular for this problem.

\begin{theorem}
\label{gendir} Let $X=\{X_{1},...,X_{m}\}$ be a system of smooth vector fields
satisfying \eqref{1.1x}. Let $\Omega_{T}\subset\mathbf{R}^{n+1}$ be an $X$-NTA
domain and let $f\in C(\partial_{p}\Omega_{T})$ and $g\in C^{0,\beta}%
(\Omega_{0})$ for some $0<\beta\leq\alpha$ and some neighborhood $\Omega
_{0}\supset\bar{\Omega}_{T}$. Then there exists a solution $u\in C^{2,\beta
}(\Omega_{T})\cap C(\partial_{p}\Omega_{T})$ to the problem in \eqref{dirp}.
\end{theorem}

\noindent\textbf{Proof:} This follows from Theorem 4.1 in Uguzzoni \cite{U}.
We only need to prove that $X$-NTA domains satisfy the exterior d-cone
criterion, which is simple, but we write it down for the reader's convenience.
By definition, this criterion is satisfied if there exists $\theta>0$ such
that for every $r<2r_{0}$ we have
\[
|B_{d}^{t_{0}-r^{2}}(x_{0},r)\backslash\cap\Omega_{T}^{C}|\geq\theta
|B_{d}(x_{0},r)|.
\]
To see that this is so, observe that if $\Omega_{T}$ is an $X$-NTA domain and
if $\underline{A}_{r}^{-}(x_{0},t_{0}) = (x_{1},t_{1})$, we have $B_{d}%
^{t_{1}}(x_{1},r/M)\subset\Omega_{T}^{C}$, and since $\underline{A}_{r}%
^{-}(x_{0},t_{0})\in B_{d}^{t_{0}-r^{2}}(x_{0},r)$, we can find a ball
$B_{d}^{t_{0}-r^{2}}(x^{\prime},r/2M)\subset B_{d}^{t_{0}-r^{2}} (x_{0}%
,r)\cap\Omega_{T}^{C}$. We now have
\[
|B_{d}^{t_{0}-r^{2}}(x^{\prime},r/2M)|\geq(4M)^{-Q}|B_{d}^{t_{0} -r^{2}%
}(x^{\prime},2r)|\geq(4M)^{-Q}|B_{d}(x_{0},r)|.
\]
\hfill\hfill\hfill$\Box$ \newline

Note that this theorem is actually proved for operators in non-divergence form
in \cite{U}, whereas our operators are in divergence form. However, since the
coefficients $a_{ij}$ are smooth, and due to the shape of the adjoints
$X_{i}^{\ast},$ see below (\ref{1.1x}), the results extend to our situation.

Now that we know that all points in $\partial_{p}\Omega_{T}$ are regular for
the Dirichlet problem in (\ref{dirp}), we recall the $X$-parabolic measure
$\omega_{X}$ from the introduction (see (\ref{caloricM}) and below).

By using the simple geometrical argument in Lemma 6.4 in \cite{LU}, one can
show that $\mathbb{R}^{n+1} \setminus C_{R}(x_{0},t_{0})$ satisfies condition
(ii) in Definition \ref{defNTA}, and thus it also satisfies the uniform
exterior $d$-cone condition, and one can solve the Dirichlet problem there
(but we stress that the set is not in general $X$-NTA). The same is true of
the intersection of two sets that satisfy condition (ii) in Definition
\ref{defNTA}. This is used to prove the following lemma (Theorem 6.5 in
\cite{LU}):

\begin{lemma}
\label{regular} Let $D \in\mathbb{R}^{n}$ be open and bounded. Then, for every
$\delta>0$ there exists a set $D_{\delta}$ such that $\{x \in D:
d(x,\partial_{p} D) > \delta\} \subset D_{\delta}\subset D$, and the cylinder
$D_{\delta}\times(t_{1},t_{2})$ satisfies the uniform exterior $d$-cone condition.
\end{lemma}

Note that the cylinder $D_{\delta}\times(t_{1},t_{2})$ is thus regular for the
Dirichlet problem, and has an $X$-parabolic measure.

\subsection{Preliminary estimates}

We will need the following results for future use. First, we have the interior
Harnack inequality, see \cite[Theorem 19.2]{BBLU} or \cite[Theorem 1.2]{CCR}:

\begin{lemma}
\label{lem4.1}\emph{(Harnack inequality)} Let $0<h_{1}<h_{2}<1$ and $\gamma
\in(0,1)$. Assume that $u\in\Gamma^{2}(B_{d}(x_{0},r)\times(t_{0}-r^{2}%
,t_{0}))\cap C(\overline{B_{d}(x_{0},r)}\times\lbrack t_{0}-r^{2},t_{0}])$.
Then there exists a positive constant $c=c(h_{1},h_{2},\gamma,r_{0},X)$ such
that for every $r\leq r_{0}$, $(x_{0},t_{0})\in\mathbb{R}^{n+1},$
\[
\max_{\overline{B_{d}(x_{0},\gamma r)}\times\lbrack t_{0}-h_{2}r^{2}%
,t_{0}-h_{1}r^{2}]}u\leq cu(x_{0},t_{0})
\]
if $u$ is a non-negative solution to the heat equation $Hu=0$ in $B_{d}%
(x_{0},r)\times(t_{0}-r^{2},t_{0})$, while%
\[
\max_{\overline{B_{d}(x_{0},\gamma r)}\times\lbrack t_{0}+h_{1}r^{2}%
,t_{0}+h_{2}r^{2}]}u\leq cu(x_{0},t_{0}),
\]
if $u$ is a non-negative solution to the adjoint heat equation $\hat{H}u=0$ in
$B_{d}(x_{0},r)\times(t_{0}-r^{2},t_{0}).$
\end{lemma}

We will also need a Cacciopoli (energy) estimate.

\begin{lemma}
\label{lem4.2b} Take a non-negative function $u$ in $C_{2r}(x,t)$ such that
$Hu\leq0$. Then we have
\[
\int_{C_{r}(x,t)}|Xu|^{2}dyds+\int_{C_{r}(x,t)}\frac{\partial(u^{2})}{\partial
s}dyds\leq\frac{c}{r^{2}}\int_{C_{2r}(x,t)}u^{2}dyds.
\]

\end{lemma}

\noindent\textbf{Proof:} By Theorem 1.5 in \cite{GN} we can find a cut-off
function $\psi\in C_{0}^{\infty}(C_{2r}(x,t))$ such that $\phi=1$ on
$C_{2r}(x,t)$ and $|X\psi|\leq c/r$. The lemma is then proved by standard
methods using partial integration, see also Lemma 3.1 in \cite{DZN}.
\hfill$\Box$ \newline

Next, we need estimates for the fundamental solution, and a solution of the
Cauchy problem in $\mathbb{R}^{n}\times\lbrack0,T]$. Theorem 10.7 in
\cite{BBLU} give, respectively,

\begin{lemma}
\label{lem4.3} There exists a fundamental solution $\Gamma(x,t,\xi,\tau)$ for
$H$ on $\mathbb{R}^{n+1}$ and constants $c_{1}$ and $c_{2}$ depending on $X$
and $T$, such that
\[
\frac{\exp(-c_{2}d(x,\xi)^{2}/(t-\tau))}{c_{1}|B_{d}(x,\sqrt{t-\tau})|}%
\leq\Gamma(x,t,\xi,\tau)\leq\frac{c_{1}\exp(-d(x,\xi)^{2}/c_{2}(t-\tau
))}{|B_{d}(x,\sqrt{t-\tau})|}%
\]
if $0<t-\tau<T$ and $x,\xi\in\Omega_{T}$. Further, $\Gamma(x,t,\xi,\tau)=0$
for $\tau-\tau\leq0$.
\end{lemma}

\begin{lemma}
\label{lem4.4} Let $\mu> 0$ and $T>0$ be such that $\mu T$ is small enough.
Take $g \in C(\mathbb{R}^{n})$ such that $|g(x)| \leq c \exp{\mu d(x,0)^{2}}$
for some constant $c>0$. Then
\[
u(x,t) = \int_{\mathbb{R}^{n}} \Gamma(x,t,\xi,0) g(\xi) d\xi
\]
is a solution to $Hu = 0$ in $\mathbb{R}^{n} \times[0,T]$ with the initial
condition $u(x,0) = g(x)$.
\end{lemma}

\section{A backward in time Harnack inequality at the boundary\label{wrangle}}

In this section we consider non-negative solutions to the subelliptic heat
equation $Hu=0$ and to the adjoint subelliptic heat equation $\hat{H}u=0$ in
$\Omega_{T}$. We note that the adjoint subelliptic heat equation is given by%
\[
\hat{H}u=\partial_{t}+\sum_{i,j=1}^{m}X_{j}^{\ast}(a_{ij}(x,t)X_{i}u).
\]
The following will be used as assumptions in most lemmas and theorems from now
on:
\begin{align}
&  \text{Let $X=\{X_{1},...,X_{m}\}$ be a system of smooth vector fields
satisfying \eqref{1.1x}.}\label{ass}\\
&  \text{Let $\Omega_{T}\subset\mathbb{R}^{n+1}$ be an $X$-NTA-domain.}%
\nonumber
\end{align}

\begin{lemma}
\label{lemmaacchar} Assume (\ref{ass}) and let $(y_{1},s_{1})$, $(y_{2}%
,s_{2})\in\Omega_{T}$; suppose that $(s_{2}-s_{1})^{1/2}\geq\eta^{-1}%
d_{p}(y_{1},s_{1},y_{2},s_{2})$ for some $\eta>1$; and that $d(y_{1}%
,\partial\Omega_{T})>\epsilon$, $d(y_{2},\partial\Omega_{T})>\epsilon$,
$(T-s_{2})>\epsilon^{2}$, $s_{1}>\epsilon^{2}$ and $d_{p}((y_{1},s_{1}%
),(y_{2},s_{2}))<c\epsilon$ for some $\epsilon>0$. Then there exists a
constant $\hat{c}=\hat{c}(H,\eta,c,r_{0})$, $\hat{c}\leq1$, such that
\[
u(y_{1},s_{1})\leq\hat{c}u(y_{2},s_{2}),
\]
if $u$ is a nonnegative solution to the heat equation $Hu=0$ in $\Omega_{T}$,
while%
\[
u(y_{2},s_{2})\leq\hat{c}u(y_{1},s_{1}),
\]
if $u$ is a nonnegative solution to the adjoint heat equation $\hat{H}u=0$ in
$\Omega_{T}.$
\end{lemma}

\noindent\textbf{Proof:} This is proved by using the parabolic Harnack chain
from Definition \ref{defNTA}, and applying Lemma \ref{lem4.1} in each
cylinder. We omit the details.\hfill$\Box$ \newline

\begin{lemma}
\label{lem4.7} Assume (\ref{ass}), let $(x_{0},t_{0})\in S_{T}$, and let
$r<\min(r_{0}/2,\sqrt{T-t_{0}}/4,\sqrt{t_{0}}/4\}$. Let $u$ be a non-negative
solution to either the subelliptic heat equation, or the adjoint subelliptic
heat equation, in $\Omega_{T}\cap C_{4r}(x_{0},t_{0})$, and assume that $u$
vanishes continuously on $\Delta(x_{0},t_{0},2r)$. Then there exist constants
$c=c(n,M,r_{0})$, $c\geq1$, and $\gamma=\gamma(H,M,r_{0})>0$, such that
\[
u(x,t)d_{p}(x,t,S_{T})^{\gamma}\leq cr^{\gamma}u(A_{r}^{+}(x_{0},t_{0})),
\]
if $u$ solves the heat equation, and%
\[
u(x,t)d_{p}(x,t,S_{T})^{\gamma}\leq cr^{\gamma}u(A_{r}^{-}(x_{0},t_{0})),
\]
if $u$ solves the adjoint heat equation, whenever $(x,t)\in\Omega_{T}\cap
C_{r}(x_{0},t_{0})$.
\end{lemma}

For a proof, see Lemma 3.1 in \cite{FGGMN}. In fact, we can use that proof
almost word for word, but be aware that $A_{r}^{+}(x_{0},t_{0})$ and
$A_{r}^{-}(x_{0},t_{0})$ are defined differently, since the domain in
\cite{FGGMN} is not time-dependent. Using a similar argument we can prove the
following lemma;

\begin{lemma}
\label{lem4.8} Assume (\ref{ass}), let $(x_{0},t_{0})\in S_{T}$, and let
$r<\min(r_{0}/2,\sqrt{T-t_{0}}/4,\sqrt{t_{0}}/4\}$. Let $u$ be a non-negative
solution to the subelliptic heat equation in $\Omega_{T}\cap C_{4r}%
(x_{0},t_{0})$, and assume that $u$ vanishes continuously on $\Delta
(x_{0},t_{0},2r)$. Then there exist constants $c=c(H,M,r_{0})$, $c\geq1$, and
$\gamma=\gamma(n,M,r_{0})>0$, such that
\[
u(A_{r}^{-}(x_{0},t_{0}))\leq cr^{\gamma}u(x,t)d_{p}(x,t,\partial_{p}%
\Omega_{T})_{{}}^{-\gamma}%
\]
whenever $(x,t)\in\Omega_{T}\cap C_{r}(x_{0},t_{0})$. If, on the other hand,
$u$ is a solution the adjoint heat equation, then, under the above
assumptions,
\[
u(A_{r}^{+}(x_{0},t_{0}))\leq cr^{\gamma}u(x,t)d_{p}(x,t,\partial_{p}%
\Omega_{T})^{-\gamma}%
\]
whenever $(x,t)\in\Omega_{T}\cap C_{r}(x_{0},t_{0})$.
\end{lemma}

\begin{lemma}
\label{doftticka} Assume (\ref{ass}). Let $(x_{0},t_{0})\in S_{T}$ and let
$r<\min\{r_{0}/2,\sqrt{T-t_{0}}/4,\sqrt{t_{0}}/4\}$. Then there exists a
constant $c=c(H,M,r_{0})$ such that
\begin{equation}
\omega_{X}(x,t,\Delta(x_{0},t_{0},2r))>c \label{prop1}%
\end{equation}
and%
\begin{equation}
\hat{\omega}_{X}(x,t,\Delta(x_{0},t_{0},2r))>c \label{prop2}%
\end{equation}
for $(x,t)\in\Omega_{T}\cap C_{r}(x_{0},t_{0})$.
\end{lemma}

\noindent\textbf{Proof:} By Lemma \ref{regular}, we can choose a set $U$ which
is regular for the Dirichlet problem and such that $B_{d}(x_{0},r) \subset U
\subset B_{d}(x_{0},2r)$. Let $C:=U\times[t_{0} - 2r^{2}, t_{0}+2r^{2}]$. Let
$y_{0}$ be defined by $(y_{0},t_{0}-2r^{2}) = \underline{A}_{r}^{-}%
(x_{0},t_{0})$. By Lemma \ref{regular} we can also find a set $U^{\prime}$
which is regular for the Dirichlet problem such that $B_{d}(y_{0},r/2M)
\subset U^{\prime}\subset B_{d}(y_{0},r/M)$. Define $C^{\prime}:= U^{\prime
}\times[t_{0} - 2r^{2}, t_{0}-2r^{2}+r^{2}/(4M^{2})]$ and $B: = U^{\prime
}\times\{t=t_{0}-2r^{2}\}$, or the bottom of $C^{\prime}$; then $C^{\prime}
\subset C \setminus\Omega_{T}$.

Now, we let $v(x,t)=\omega_{X}(x,t,B,C)$ and $v^{\prime}(x,t)=\omega
_{X}(x,t,B,C^{\prime})$. By the maximum principle, we have $\omega
_{X}(x,t,\Delta(x_{0},t_{0},2r)\geq v(x,t)$ in $\Omega_{T}\cap U$ and thus
also in $\Omega_{T}\cap C_{r}(x_{0},t_{0})$, and $v(x,t)\geq v^{\prime}(x,t)$
in $C^{\prime}$. By the Harnack principle, we have
\[
\inf_{\Omega_{T}\cap C_{r}(x_{0},t_{0})}\omega_{X}(x,t,\Delta(x_{0}%
,t_{0},2r)\geq\inf_{\Omega_{T}\cap C_{r}(x_{0},t_{0})}v(x,t)\geq
\]%
\begin{equation}
c^{-1}v\left(  y_{0},t_{0}-2r^{2}+\frac{r^{2}}{8M^{2}}\right)  \geq
c^{-1}v^{\prime}\left(  y_{0},t_{0}-2r^{2}+\frac{r^{2}}{8M^{2}}\right)  .
\label{run}%
\end{equation}
We can extend the function $v^{\prime}$ to the cylinder
\[
\tilde{C}=U^{\prime}\times\left[  t_{0}-2r^{2}-\frac{r^{2}}{4M^{2}}%
,t_{0}-2r^{2}+\frac{r^{2}}{4M^{2}}\right]
\]
by setting
\[
\tilde{v}(x,t)=\omega_{X}(x,t,\partial_{p}\tilde{C}\cap\{t\leq t_{0}%
-2r^{2}\},\tilde{C}),
\]
that is, extending $v^{\prime}$ by setting it as $1$ below $B$. We now apply
the Harnack inequality to $v^{\prime}$ in $\tilde{C}$ and obtain
\begin{equation}
v^{\prime}\left(  y_{0},t_{0}-2r^{2}+\frac{r^{2}}{8M^{2}}\right)  =\tilde
{v}\left(  y_{0},t_{0}-2r^{2}+\frac{r^{2}}{8M^{2}}\right)  \geq c^{-1}%
\tilde{v}(y_{0},t_{0}-2r^{2})=c^{-1}. \label{run2}%
\end{equation}
Combining (\ref{run}) and (\ref{run2}) and (\ref{prop1}) follows.

To see that (\ref{prop2}) also holds, let $y_{0}$ be such that $(y_{0}%
,t_{0}+2r^{2})=\underline{A}_{r}^{+}(x_{0},t_{0})$ and let $U^{\prime}$ be a
set which is regular for the Dirichlet problem, such that $B_{d}%
(y_{0},r/2M)\subset U^{\prime}\subset B_{d}(y_{0},r/M)$. Then we define
$C^{\prime}:=U^{\prime}\times\lbrack t_{0}+2r^{2}-r^{2}/(4M^{2}),t_{0}%
+2r^{2}]$ and $B:=U^{\prime}\times\{t=t_{0}+2r^{2}\},$ the top of $C^{\prime
}.$ Now we can argue in line with the proof of (\ref{prop1}). We omit the
details. \hfill$\Box$ \newline

\begin{remark}
\label{kropotkin} Note that we can actually extend Lemma \ref{doftticka} to
the set $\{(x,t) \in\partial\Omega_{T}: t=0\}$. The proof there is much
simpler, but uses the same idea.
\end{remark}

\begin{lemma}
\label{lem4.5} Assume (\ref{ass}). Let $(x_{0},t_{0})\in S_{T}$ and let
$r<\min\{r_{0}/2,\sqrt{T-t_{0}}/4,\sqrt{t_{0}}/4\}$. Let $u$ be a non-negative
solution to either the subelliptic heat equation or the adjoint subelliptic
heat equation in $\Omega_{T}\cap C_{4r}(x_{0},t_{0})$. Assume that $u$
vanishes continuously on $\partial_{p}\Omega_{T}\cap C_{2r}(x_{0},t_{0})$.
Then there exist constants $c=c(H,M,r_{0})$, $c\geq1$, and $\alpha
=\alpha(H,M,r_{0})\in(0,1)$, such that
\[
u(x,t)\leq c\biggl (\frac{d_{p}(x,t,x_{0},t_{0})}{r}\biggr )^{\alpha}%
\sup_{\Omega_{T}\cap C_{2r}(x_{0},t_{0})}u
\]
whenever $(x,t)\in\Omega_{T}\cap C_{r}(x_{0},t_{0})$.
\end{lemma}

\noindent\textbf{Proof:} First, assume that $u$ is a solution to the
subelliptic heat equation. We will prove that there exists a $\Theta
=\Theta(H,M,r_{0})\in(0,1)$ such that
\[
\sup_{\Omega_{T}\cap C_{r}(x_{0},t_{0})}u\leq\Theta\sup_{\Omega_{T}\cap
C_{2r}(x_{0},t_{0})}u.
\]
By Lemma \ref{regular}, we can choose a set $U$ which is regular for the
Dirichlet problem and such that $B_{d}(x_{0},3r/2)\subset U\subset B_{d}%
(x_{0},2r)$, and let $C=U\times\lbrack t_{0}-4r^{2},t_{0}+4r^{2}]$. If
$(x,t)\in\Omega_{T}\cap C_{r}(x_{0},t_{0})$ we have%

\begin{align}
u(x,t)  &  =\int_{\partial_{p}(\Omega_{T}\cap C)}ud\omega(x,t,y,s,\Omega
_{T}\cap C)=\int_{\Omega_{T}\cap\partial_{p}C}ud\omega(x,t,y,s,\Omega_{T}\cap
C)\leq\nonumber\\
\text{ \ }  &  \leq\left(  \sup_{\Omega_{T}\cap C_{2r}(x_{0},t_{0})}u\right)
\omega(x,t,\Omega_{T}\cap\partial_{p}C,\Omega_{T}\cap C)=\nonumber\\
&  =\left(  \sup_{\Omega_{T}\cap C_{2r}(x_{0},t_{0})}u\right)  (1-\omega
(x,t,\partial_{p}\Omega_{T}\cap C,\Omega_{T}\cap C)). \label{bulb}%
\end{align}
It is enough to prove that $\omega(x,t,\partial_{p}\Omega_{T}\cap C,\Omega
_{T}\cap C)>c$, and this is done in the same way as in the proof of Lemma
\ref{doftticka}.

We now iterate the procedure, and get $u(x,t)\leq\Theta^{2}\sup_{\Omega
_{T}\cap C_{2r}(x_{0},t_{0})}u$ if $(x,t)\in\Omega_{T}\cap C_{r/2}^{-}%
(x_{0},t_{0})$. Continuing like that, $u(x,t)\leq\Theta^{k}\sup_{\Omega
_{T}\cap C_{2r}(x_{0},t_{0})}u$ if $(x,t)\in\Omega_{T}\cap C_{2^{-k+1}r}%
(x_{0},t_{0})$. In fact, $u(x,t)\leq\Theta^{k}\sup_{\Omega_{T}\cap
C_{2r}(x_{0},t_{0})}u$ if $2^{-k}r\leq d_{p}(x,t,x_{0},t_{0})\leq2^{-k+1}r$.
We want to show that
\[
\Theta^{k}\leq\left(  \frac{d_{p}(x,t,x_{0},t_{0})}{r}\right)  ^{\alpha}%
\]
when $2^{-k}r\leq d_{p}(x,t,x_{0},t_{0})\leq2^{-k+1}r$. In fact, we have
$(d_{p}(x,t,x_{0},t_{0})/r)^{\alpha}\geq2^{-\alpha k}>(2^{-\alpha})^{k}$, so
we let $\alpha$ satisfy $2^{-\alpha}=\Theta$, which completes the proof when
$u$ is a solution to the subelliptic heat equation. If $u$ is a solution to
the adjoint subelliptic heat equation, the proof is similar, the only
difference is that we must replace $\omega$ in (\ref{bulb}) by $\hat{\omega}%
$.\hfill$\Box$ \newline

\begin{lemma}
\label{lem4.6} Assume (\ref{ass}). Let $(x_{0},t_{0})\in S_{T}$ and let
$r<\min\{r_{0}/2,\sqrt{T-t_{0}}/4,\sqrt{t_{0}}/4\}$. Let $u$ be a non-negative
solution to either the subelliptic heat equation or the adjoint subelliptic
heat equation in $\Omega_{T}\cap C_{4r}(x_{0},t_{0})$. Assume that $u$
vanishes continuously on $\Delta(x_{0},t_{0},2r)$. Then, whenever
$(x,t)\in\Omega_{T}\cap C_{r}(x_{0},t_{0})$, there exists a constant
$c=c(H,M,r_{0})$, $c \geq1$, such that
\[
u(x,t)\leq cu(P_{r}(x_{0},t_{0})),
\]
where $P_{r}(x_{0},t_{0})=A_{r}^{+}(x_{0},t_{0})$ if $u$ solves the
subelliptic heat equation and $P_{r}(x_{0},t_{0})=A_{r}^{-}(x_{0},t_{0})$ if
$u$ solves the adjoint subelliptic heat equation.
\end{lemma}

\noindent\textbf{Proof:} This follows from Lemma \ref{lem4.1}, Lemma
\ref{lem4.5} and a classical argument by contradiction, see the proof of
Theorem 3.1 in \cite{S}.\hfill$\Box$ \newline

\begin{lemma}
\label{smultron} Assume (\ref{ass}). Let $(x_{0},t_{0})\in S_{T}$ and let
$r<\min\{r_{0}/2,\sqrt{T-t_{0}}/4,\sqrt{t_{0}}/4\}$. Let $u$ be a non-negative
solution to the subelliptic heat equation in $\Omega_{T}$. Assume that $u$
vanishes continuously on $\partial_{p} \Omega_{T} \setminus\Delta(x_{0}%
,t_{0},r/2)$. Then there exists a constant $c=c(H,M,r_{0})$ such that
\[
u(x,t)\leq c u(A^{+}_{r}(x_{0},t_{0}))
\]
whenever $(x,t)\in\Omega_{T} \setminus C_{r}(x_{0},t_{0})$.
\end{lemma}

\noindent\textbf{Proof:} On account of the maximum principle, it is enough to
prove the statement on the set $D = \partial C_{r}(x_{0},t_{0}) \cap\Omega
_{T}$. This set can be divided into the two sets $A = \{(x,t) \in D:
d_{p}(x,t,\partial_{p} \Omega_{T}) < r/8 \}$ and $B = D \setminus A$. On the
set $B $, we use the interior Harnack principle. If $(x,t) \in A$, we have
$(x,t) \in C_{r/8}(P)$, where $P \in\partial_{p} \Omega_{T}$ is a point such
that $d_{p}(x,t,P) = d_{p}(x,t,\partial_{p} \Omega_{T})$. Since $u$ vanishes
continuously on $\Delta(P,r/4)$, we can use Lemma \ref{lem4.6} to draw the
conclusion that $u(x,t) \leq c u(A^{+}_{r/8}(P))$. We know that $d_{p}%
(A^{+}_{r/8}(P),\partial_{p} \Omega_{T}) > r/8M$, so we can use the interior
Harnack principle again to conclude that $u(A^{+}_{r/8}(P)) \leq c u(A^{+}%
_{r}(x_{0},t_{0}))$. \hfill$\Box$ \newline

\begin{lemma}
\label{tag} Assume (\ref{ass}). Let $(x_{0},t_{0})\in S_{T}$ and let
$r<\min\{r_{0}/2,\sqrt{T-t_{0}}/4,\sqrt{t_{0}}/4\}$. Let $u$ be a non-negative
solution to the subelliptic heat equation in $\Omega_{T}$. Assume that $u$
vanishes continuously on $\partial_{p} \Omega_{T} \setminus\Delta(x_{0}%
,t_{0},r/2)$. Then there exists a constant $c = c(H,M,r_{0})$ such that
\[
u(x,t)\leq c u(A^{+}_{r}(x_{0},t_{0})) \omega_{X} (x,t,\Delta(x_{0}%
,t_{0},2r))
\]
whenever $(x,t)\in\Omega_{T} \setminus C_{r}(x_{0},t_{0})$.
\end{lemma}

\noindent\textbf{Proof:} This follows from Lemma \ref{smultron} and Lemma
\ref{doftticka}. \hfill$\Box$ \newline

\begin{lemma}
\label{lem4.9} Assume (\ref{ass}). Let $u$ be a solution to either the
subelliptic heat equation or to the adjoint subelliptic heat equation in
$\Omega_{T}$. Assume that $u$ vanishes continuously on $S_{T}$. Let
$\delta<4r_{0}$ and let $D_{\delta}$ denote the closure of the set
$\{(x,t)\in\Omega_{T}:d_{p}(x,t,\partial_{p}\Omega_{T})>\delta\}$; then there
exists a constant $c=c(H,M,r_{0},\delta,\text{diam}(\Omega_{T}),T)$ such that
\[
\max_{D_{\delta}}u(x,t)\leq c\min_{D_{\delta}}u(x,t).
\]

\end{lemma}

\noindent\textbf{Proof:} First, let $u$ be a solution to the subelliptic heat
equation. Since $u$ is continuous on the compact set $D_{\delta}$, there exist
$(x_{1},t_{1})$ and $(x_{2},t_{2})$ in $D_{\delta}$ such that $\min
_{D_{\delta}}u(x,t)=u(x_{1},t_{1})$ and $\max_{D_{\delta}}u(x,t)=u(x_{2}%
,t_{2})$.

We define $\Omega_{\lbrack\delta/2,T]}=\{(x,t)\in\Omega:t\in((\delta
/2)^{2},T)\}$, $A=\{(x,t)\in\Omega:t=(\delta/2)^{2}\}$ and
\[
B=\left\{  (x,t)\in\Omega:t=\left(  \frac{\delta}{2}\right)  ^{2}+2\left(
\frac{\delta}{8}\right)  ^{2},d_{p}(x,t,\partial\Omega)>\frac{\delta}%
{8M}\right\}  .
\]
We have
\[
u(x_{2},t_{2})\leq\sup_{\Omega_{\lbrack\delta/2,T]}}u\leq\sup_{A}u
\]
by the maximum principle, and we aim to prove that $\sup_{A}u\leq c\sup_{B}u$.
If $(x,t)\in A$ and $d_{p}(x,t,\partial\Omega_{T})\geq\delta/8$, then we can
apply the Harnack principle to reach a point in $B$. If $(x,t)$ is closer to
the boundary, we can apply Lemma \ref{lem4.6} with $r=\delta/8$. Finally, we
take a point $(y,s)\in\bar{B}$ such that $u(y,s)=\max_{\bar{B}}u $, and by the
Harnack principle we have $u(y,s)\leq cu(x_{1},t_{1})$, which completes the
proof for solutions to the subelliptic heat equation. If $u$ is a solution to
the adjoint subelliptic heat equation, we replace $\Omega_{\lbrack\delta
/2,T]}$ with $\Omega_{\lbrack0,T-\delta/2]}$, and let $A=\{(x,t)\in
\Omega:t=T-(\delta/2)^{2}\}$ and
\[
B=\left\{  (x,t)\in\Omega:t=T-\left(  \frac{\delta}{2}\right)  ^{2}-2\left(
\frac{\delta}{8}\right)  ^{2},d_{p}(x,t,\partial\Omega)>\frac{\delta}%
{8M}\right\}  .
\]
Then we argue as above, and this completes the proof.\hfill$\Box$ \newline

Now we prove a technical lemma that we need in order to prove the backward
Harnack inequality. Let
\[
\Omega_{1}=C_{Kr,r}(x_{0},t_{0})\cap\Omega_{T}.
\]
By Lemma \ref{regular}, we can find a set $\tilde{B}$ which is regular for the
Dirichlet problem and such that
\[
B_{d}(x_{0},(2K-1)r)\subset\tilde{B}\subset B_{d}(x_{0},2Kr).
\]
Further, we set $\Omega_{2}=[\tilde{B}\times(t_{0}-3r^{2},t_{0}+3r^{2}%
)]\cap\Omega_{T}$.

\begin{lemma}
\label{choklad} Assume (\ref{ass}). Let $u$ be a non-negative solution to
either the subelliptic heat equation or the adjoint subelliptic heat equation
in $\Omega_{T}$. Assume that $u$ vanishes continuously on $S_{T}$. Let
$0<\delta\ll\sqrt{T}$ be a fixed constant, and let $(x_{0},t_{0})\in S_{T}$,
$\delta^{2}\leq t_{0}\leq T-\delta^{2}$. Then there exists ${K}={K}(H,M)$ such
that if
\[
\sup_{\Omega_{1}}u\geq(2{K})^{-\gamma}\sup_{\Omega_{2}}u
\]
then
\[
\sup_{\bar{\Omega}_{2}\cap\{t=t_{0}-3r^{2}\}}u\geq\frac{1}{2}\sup_{\Omega_{1}%
}u.
\]

\end{lemma}

\noindent\textbf{Proof:} Assume that $u$ is a solution to the subelliptic heat
equation. For simplicity of writing, we assume that $\sup_{\Omega_{2}}u=1$,
and that $t_{0}-3r^{2}=0$. We also set $\sup_{\Omega_{1}}u=M_{0}$. Our
assumption now says that $M_{0}\geq(2{K})^{-\gamma}$, and we will prove that
$\sup_{\bar{\Omega}_{2}\cap\{t=0\}}u\geq M_{0}/2$ if ${K}$ is large enough. We
will prove this by contradiction, so we assume that $\sup_{\bar{\Omega}%
_{2}\cap\{t=0\}}u<M_{0}/2$.

By using the $X$-parabolic measure we can write
\[
u(x,t)=\int_{\partial_{p}\Omega_{2}}u(y,s)d\omega(x,t,y,s,\Omega_{2}).
\]
Now, define $\Gamma_{1}=\bar{\Omega}_{2}\cap\{t=0\}$ and $\Gamma_{2}%
=\partial_{p}\Omega_{2}\setminus\{\Gamma_{1}\cup S_{T}\}$. If $(x,t)\in
\Omega_{2}$, using the maximum principle, we have
\begin{equation}
u(x,t)\leq(\sup_{\Gamma_{1}}u)\omega(x,t,\Gamma_{1},\Omega_{2})+(\sup
_{\Gamma_{2}}u)\omega(x,t,\Gamma_{2},\Omega_{2})\leq\newline\frac{M_{0}}%
{2}\omega(x,t,\Gamma_{1},\Omega_{2})+\omega(x,t,\Gamma_{2},\Omega_{2}).
\label{beb}%
\end{equation}
Let $\phi\in C_{0}^{\infty}(\mathbb{R}^{n})$ be such that $0\leq\phi\leq1$,
$\phi=1$ on $B_{d}(x_{0},2({K}+2)r)\setminus B_{d}(x_{0},2({K}-2)r)$ and
$\phi=0$ on $B_{d}(x_{0},2({K}-3)r)\cup(\mathbb{R}^{n}\setminus B_{d}%
(x_{0},2({K}+3)r))$. Let $\Phi$ be the solution to the Cauchy problem with
$\phi$ as initial data, i.e., recalling Lemma \ref{lem4.4} we set
\[
\Phi(x,t)=\int_{\mathbb{R}^{n}}\Gamma(x,t,\xi,0)\phi(\xi)d\xi,
\]
for $x\in\mathbb{R}^{n}$, $t\geq0$.

As in Lemma \ref{doftticka}, we prove that $\Phi(x,t)\geq c$ for
$(x,t)\in\Gamma_{2}$.
By the maximum principle, Lemma \ref{lem4.4} and the fact that $\omega$ is a
probability measure, we then have $\Phi(x,t)\geq c\omega(x,t,\Gamma_{2}%
,\Omega_{2}) $ on $\Omega_{2}$.

If $(x,t)\in\bar{\Omega}_{1}$, then by Lemma \ref{lem4.3}
\[
\omega(x,t,\Gamma_{2},\Omega_{2})\leq c^{-1}\Phi(x,t)\leq c^{-1}\int%
_{B_{d}(x_{0},2({K}+3)r)\setminus B_{d}(x_{0},2({K}-3)r)}\exp\left(
-cd(x,y)^{2}/r^{2}\right)  dy.
\]
Since $x\in B_{d}(x_{0},{K}r)$ and $y\in B_{d}(x_{0},2({K}+3)r)\setminus
B_{d}(x_{0},2({K}-3)r)$, we have
\[
d(x,y)\geq d(x_{0},y)-d(x_{0},x)\geq2({K}-3)r-{K}r=r({K}-6)\geq\tilde{C}{K}r
\]
if ${K}$ is large enough. Hence,
\[
\omega(x,t,\Gamma_{2},\Omega_{2})\leq c^{-1}\int_{\{y:d(x,y)\geq\tilde{C}%
{K}r\}}\exp\left(  -cd(x,y)^{2}/r^{2}\right)  dy.
\]
By means of the estimate $C|x-y|\leq d(x,y)$ and integration on the level sets
of $|x-y|$, one can prove that
\begin{equation}
\omega(x,t,\Gamma_{2},\Omega_{2})\leq\exp(-c{K}^{2}), \label{beb2}%
\end{equation}
if ${K}$ is large enough.

Let $(x_{1},t_{1})\in\Omega_{1}$ be a point such that $u(x_{1},t_{1})=M_{0} $.
Combining (\ref{beb}) and (\ref{beb2}) we get the following estimate for
$M_{0}$ if ${K}$ is large enough;
\[
M_{0}=u(x_{1},t_{1})\leq M_{0}/2+c^{-1}\exp(-c{K}^{2})<M_{0}/2+(2
{K})^{-\gamma}/2.
\]
That is, $M_{0}<(2{K})^{-\gamma}$, which is a contradiction.

When $u$ is a solution to the adjoint subelliptic heat equation, we must use
$\hat{\omega}$ instead of $\omega.$ We simplify notation as above, with the
only difference is that we set $t_{0}+3r^{2}=T.$ Then we define $\Gamma
_{1}=\bar{\Omega}_{2}\cap\{t=T\}$ and $\Gamma_{2}=\partial_{p}\Omega
_{2}\setminus\{\Gamma_{1}\cup S_{T}\}$ and argue as before, using that
$\Gamma^{\hat{H}}(x,t,\xi,\tau)=\Gamma^{H}(\xi,\tau,x,t)$ is a fundamental
solution for the adjoint problem. This completes the proof.\hfill$\Box$

We can now prove our backward Harnack inequality:

\begin{theorem}
\label{lem4.12} Assume (\ref{ass}). Let $u$ be a non-negative solution to
either the subelliptic heat equation or the adjoint subelliptic heat equation
in $\Omega_{T}$. Assume that $u$ vanishes continuously on $S_{T}$. Let
$0<\delta_{0}\ll\sqrt{T}$ be a fixed constant, and let $(x_{0},t_{0})\in
S_{T}$, $\delta_{0}^{2}\leq t_{0}\leq T-\delta_{0}^{2}$. Then there exists a
constant $c=c(H,M,r_{0},\text{diam}(\Omega_{T}),T,\delta)$, $c\geq1$, such
that
\[
u(x,t)\leq cu(A_{r}(x_{0},t_{0}))
\]
whenever $(x,t)\in\Omega_{T}\cap C_{r/4}(x_{0},t_{0})$, and for $r<\min
\{r_{0}/2,\sqrt{(T-t_{0}-\delta^{2})/4},\sqrt{(t_{0}-\delta^{2})/4}\}.$
\end{theorem}

\noindent\textbf{Proof:} Assume that $u$ is a solution to the subelliptic heat
equation. Fix $r$ as in the statement of the theorem. Define
\[
f(\alpha)=\alpha^{-\gamma}\sup_{C_{\alpha}(x_{0},t_{0})\cap\Omega_{T}}u(x,t)
\]
where $\gamma$ is the constant in Lemma \ref{lem4.8}, and let $\hat{r}%
=\max\{\alpha:r\leq\alpha\leq\delta_{0},f(\alpha)\geq f(r)\}$. By definition
of $\hat{r}$, we have
\[
\sup_{C_{r}(x_{0},t_{0})}u(x,t)\leq(r/\hat{r})^{\gamma}\sup_{C_{\hat{r}}%
(x_{0},t_{0})\cap\Omega_{T}}u(x,t).
\]
By Lemma \ref{lem4.8}, we have $u(A_{\hat{r}}^{-}(x_{0},t_{0}))\leq c\hat
{r}^{\gamma}u(A_{r}(x_{0},t_{0}))d_{p}(A_{r}(x_{0},t_{0}),\partial_{p}%
\Omega_{T})^{-\gamma}$, since $A_{r}(x_{0},t_{0})\in C_{\hat{r}}(x_{0},t_{0}%
)$. If we can prove that
\[
\sup_{C_{\hat{r}}(x_{0},t_{0})\cap\Omega_{T}}u(x,t)\leq cu(A_{\hat{r}}%
^{-}(x_{0},t_{0})),
\]
we are done, since $r/M\leq d_{p}(A_{r}(x_{0},t_{0}),\partial_{p}\Omega_{T})$,
so that $r^{\gamma}\leq M^{\gamma}d_{p}(A_{r}(x_{0},t_{0}),\partial_{p}%
\Omega_{T})^{\gamma}$. By Lemma \ref{lem4.6} we have
\[
\sup_{C_{\hat{r}}(x_{0},t_{0})\cap\Omega_{T}}u(x,t)\leq cu(A_{\hat{r}}%
^{+}(x_{0},t_{0})).
\]
Finally, if $\delta_{0}/2{K}<\hat{r}$, where ${K}$ is the constant in Lemma
\ref{choklad}, we can use Lemma \ref{lem4.9} to get $u(A_{\hat{r}}^{+}%
(x_{0},t_{0}))\leq cu(A_{\hat{r}}^{-}(x_{0},t_{0}))$, and we are finished.

On the other hand, if $r\leq\hat{r}<\delta_{0}/2{K}$, we can use Lemma
\ref{choklad}. By the definition of $\hat{r}$, we have $f(\hat{r})>f(2{K}%
\hat{r})$, that is
\[
\sup_{C_{\hat{r}}(x_{0},t_{0})\cap\Omega_{T}}u(x,t)\geq(2{K})^{-\gamma}%
\sup_{C_{2{K}\hat{r}}(x_{0},t_{0})\cap\Omega_{T}}u(x,t).
\]
Let $\Omega_{1}$ and $\Omega_{2}$ be as in Lemma \ref{choklad}, but with $r$
replaced with $\hat{r}$. We then have
\[
\sup_{\Omega_{1}}u\geq(2{K})^{-\gamma}\sup_{\Omega_{2}}u,
\]
since $C_{\hat{r}}(x_{0},t_{0})\cap\Omega_{T}\subset\Omega_{1}$ and
$\Omega_{2}\subset C_{2{K}\hat{r}}(x_{0},t_{0})\cap\Omega_{T}$. By Lemma
\ref{choklad}, this implies that there exists a point $P=(p_{x},p_{t}%
)\in\Omega_{T}$ with $p_{t}=t_{0}-3\hat{r}^{2}$ and
\[
u(P)\geq\frac{1}{2}\sup_{\Omega_{1}}u
\]
We now have
\[
\sup_{C_{\hat{r}}(x_{0},t_{0})\cap\Omega_{T}}u\leq\sup_{\Omega_{1}}%
u\leq2u(P)\leq cu(A_{\hat{r}}^{-}(x_{0},t_{0}))
\]
where the last step follows from the Harnack inequality (using Lemma
\ref{lem4.6} first if necessary). The proof is similar when $u$ is a solution
to the adjoint subelliptic heat equation, the role of $A_{r}^{+}$ and
$A_{r}^{-}$ being interchanged.\hfill$\Box$

\section{Dahlberg estimates and the doubling property for parabolic measures
\label{SecDahl}}

\noindent\textbf{Proof of Theorem \ref{lem4.10}:} We begin by proving that
\begin{equation}
c^{-1}|B_{d}(x,r)|G(x,t,A_{r}^{+}(x_{0},t_{0}))\leq\omega(x,t,\Delta
(x_{0},t_{0},r/2)), \label{feanor}%
\end{equation}
for $t\geq t_{0}$. Let $\Gamma$ be the fundamental solution to $H$, and $G$
the corresponding Green function. By definition, we have
\[
G(x,t,A_{r}^{+}(x_{0},t_{0}))=\Gamma(x,t,A_{r}^{+}(x_{0},t_{0}))-\int%
_{\partial_{p}\Omega_{T}}\Gamma(y,s,A_{r}^{+}(x_{0},t_{0}))d\omega(x,t,y,s),
\]
so that $G(x,t,A_{r}^{+}(x_{0},t_{0}))\leq\Gamma(x,t,A_{r}^{+}(x_{0},t_{0}))$.
Now, take $\alpha(M)$ such that $C_{2\alpha r}(A_{r}^{+}(x_{0},t_{0}%
))\subset\Omega_{T}$, and define
\[
A=\{(x,t)\in\Omega_{T}:t=t_{0}+2r^{2}\}\setminus C_{\alpha r/2}(A_{r}%
^{+}(x_{0},t_{0})),
\]%
\[
B=\{(x,t)\in\partial C_{\alpha r/2}(A_{r}^{+}(x_{0},t_{0})):t>t_{0}+2r^{2}\}.
\]
If $(x,t)\in A$, then $G(x,t,A_{r}^{+}(x_{0},t_{0}))=0$ by Lemma \ref{lem4.3},
so (\ref{feanor}) is true in $A$. If $(x,t)\in B$, then Lemma \ref{lem4.3} and
$(\ref{axel})$ imply that $|B_{d}(x,r)|\Gamma(x,t,A_{r}^{+}(x_{0},t_{0}))\leq
C$. We now prove that $\omega(x,t,\Delta(x_{0},t_{0},r/2))\geq C$ in $B$. By
Lemma \ref{doftticka}, we have $\omega(A_{r/4}(x_{0},t_{0}),\Delta(x_{0}%
,t_{0},r/2))\geq c^{\prime}$, and if $(x,t)\in B$ we can use the Harnack
principle, Lemma \ref{lemmaacchar}, to obtain
\[
\omega(x,t,\Delta(x_{0},t_{0},r/2))\geq c^{-1}\omega(A_{r/4}(x_{0}%
,t_{0}),\Delta(x_{0},t_{0},r/2))\geq c^{-1}c^{\prime}.
\]
So by the maximum principle, (\ref{feanor}) holds in the part of $\Omega_{T}$
bounded below by $A$ and $B$.

We now prove that $\omega(x,t,\Delta(x_{0},t_{0},r/2))\leq c|B_{d}%
(x,r)|G(x,t,A_{r}^{-}(x_{0},t_{0}))$ for $t\geq t_{0}$. Using Theorem 1.5 in
\cite{GN} we can choose $\phi\in C_{0}^{\infty}(\mathbb{R}^{n+1})$ such that
$\phi=1$ on $C_{r/2}(x_{0},t_{0})$ and $\phi=0$ on $\mathbb{R}^{n+1}\setminus
C_{3r/4}(x_{0},t_{0})$, and moreover $|\phi_{t}|\leq c/r^{2}$, $|X\phi|\leq
c/r$. Then we have
\[
\omega(x,t,\Delta(x_{0},t_{0},r/2))\leq\int_{\partial_{p}\Omega_{T}}%
\phi(y,s)d\omega(x,t,y,s).
\]
Since $|t-t_{0}|\geq5r^{2}$ (and thus $\phi(x,t)=0$) we have
\begin{multline*}
\int_{\partial_{p}\Omega_{T}}\phi(y,s)d\omega(x,t,y,s)=\\
=\int_{\Omega_{T}}XG((x,t,y,s)\cdot X\phi(y,s)+G(x,t,y,s)\phi_{s}%
(y,s)dyds=(\ast).
\end{multline*}
By the H{\"{o}}lder inequality, we have
\begin{multline*}
(\ast)\leq Cr^{n/2}\left(  \int_{C_{3r/4}(x_{0},t_{0})\cap\Omega_{T}}%
|XG|^{2}dyds\right)  ^{1/2}+\\
+Cr^{n/2-1}\left(  \int_{C_{3r/4}(x_{0},t_{0})\cap\Omega_{T}}|G|^{2}%
dyds\right)  ^{1/2}=(\ast\ast).
\end{multline*}
By using similar methods as in Lemma \ref{lem4.2b}, one can prove the
following energy estimate at the boundary: if $\alpha>1$, $u\leq0$, $Hu=0$ in
$C_{\alpha r}(x_{0},t_{0})\cap\Omega_{T}$, then
\[
\int_{C_{r}(x_{0},t_{0})\cap\Omega_{T}}|Xu|^{2}dyds\leq c/r^{2}\int_{C_{\alpha
r}(x_{0},t_{0})\cap\Omega_{T}}u^{2}dyds.
\]
Using this, we have
\[
(\ast\ast)=cr^{n/2-1}\left(  \int_{C_{r}(x_{0},t_{0})\cap\Omega_{T}%
}G(x,t,y,s)^{2}dyds\right)  ^{1/2}=(\ast\ast\ast).
\]
Given that $|t-t_{0}|\geq5r^{2}$, we know that $G(x,t,y,s)$ is a solution to
the adjoint equation $\hat{H}u(y,s)=0$ in $C_{r}(x_{0},t_{0})$, so by Lemma
\ref{lem4.6} we have $G(x,t,y,s)\leq CG(x,t,A_{r}^{-}(x_{0},t_{0}))$ for
$(y,s)\in C_{r}(x_{0},t_{0})$. This, together with (\ref{motaggsvamp}),
implies that
\[
(\ast\ast\ast)\leq C|B_{d}(x,r)|G(x,t,A_{r}^{-}(x_{0},t_{0})),
\]
and we are finished. The adjoint case, that is, bounds for $\hat{\omega},$ can
be treated treated analogously, once we have noted that $\Gamma^{\hat{H}%
}(x,t,\xi,\tau)=\Gamma^{H}(\xi,\tau,x,t)$ is a fundamental solution for the
adjoint problem. \hfill$\Box$ \newline

\begin{remark}
\label{snippet} We note that, by the same methods, we can prove the following
parts of Theorem \ref{lem4.10};
\[
\omega(x,t,\Delta(x_{0},t_{0},r/2))\leq c|B_{d}(x,r)|G(x,t,A_{r}^{-}%
(x_{0},t_{0})),
\]
and that%
\[
\hat{\omega}(x,t,\Delta(x_{0},t_{0},r/2))\leq c|B_{d}(x,r)|G(A_{r}^{+}%
(x_{0},t_{0}),x,t),
\]
if we replace the assumption on $t$ with the assumption that $(x,t)\in
\Omega_{T}$ and $x\not \in B_{d}(x_{0},4r)$.
\end{remark}

\begin{lemma}
\label{lem4.10b} Assume (\ref{ass}). Let $0<\delta<< \sqrt{T}$ be a fixed
constant, and let $(x_{0},t_{0})\in S_{T}$, $\delta^{2}\leq t_{0}\leq
T-\delta^{2}$. Let $r<\frac{1}{10} \min\{r_{0},\sqrt{T-t_{0}-\delta^{2}}%
,\sqrt{t_{0}-\delta^{2}}\}$. Assume that $(x,t)\in\Omega_{T}$ and that
$|t-t_{0}|\geq5r^{2}$. Then there exists $c=c(H,M,r_{0},\text{diam}(\Omega
_{T}),T,\delta)$, $1\leq c<\infty$, such that
\begin{align*}
&  G (x,t, A_{r}^{-}(x_{0},t_{0}) )\leq cG (x, t, A_{r}^{+}(x_{0},t_{0})
)\mbox{ if } t \geq t_{0},\\
&  G ( A_{r}^{+}(x_{0},t_{0}),x,t )\leq cG ( A_{r}^{-}(x_{0},t_{0}),x,t
)\mbox{ if }t \leq t_{0}.
\end{align*}

\end{lemma}

\noindent\textbf{Proof:} We begin with the second statement. Let $A_{r}%
^{-}(x_{0},t_{0})=(y,s)$ and $\underline{A}_{r}^{-}(x_{0},t_{0})=(y^{\prime
},s^{\prime})$. We note that since $d(y,x_{0})<r$ and $d(y^{\prime},x_{0})<r$,
there exists a path between $y$ and $y^{\prime}$ of length less than $2r$.
Since $y^{\prime}\in\Omega_{T}^{C}$, there exists a point $P\in\partial
\Omega_{T}\cap\{t=t_{0}-2r^{2}\}$ such that $d_{p}(A_{r}^{-}(x_{0}%
,t_{0}),P)<2r$. Replacing $r$ with $\sqrt{2}r$, we get a point $P^{\prime
}=(x_{1},t_{1})\in\partial\Omega_{T}\cap\{t=t_{0}-4r^{2}\}$ such that
$d_{p}(A_{\sqrt{2}r}^{-}(x_{0},t_{0}),P^{\prime})<2\sqrt{2}r$. Now let
$A_{r}^{+}(x_{0},t_{0})=(x^{\prime},t^{\prime})$. We have $d(x_{1},x^{\prime
})\leq d(x_{1},y^{\prime})+d(y^{\prime},x_{0})+d(x_{0},x^{\prime})<2r+r+r=4r$
and $t^{\prime}-t_{1}=2r^{2}+4r^{2}=6r^{2}$, so that $A_{r}^{+}(x_{0}%
,t_{0})\in C_{4r}(P)$. By Lemma \ref{lem4.12}, we have $G(A_{r}^{+}%
(x_{0},t_{0}),x,t)\leq cG(A_{4r}(P),x,t)$. We can then use the Harnack
principle to get from $A_{4r}(P)$ to $A_{r}^{-}(x_{0},t_{0})$.

To get the first statement in the lemma, we need to use Lemma \ref{lem4.12}
for the adjoint subelliptic heat equation. \hfill$\Box$ \newline

We now have all the tools we need to prove the doubling property of the
caloric measure:

\noindent\textbf{Proof of Theorem \ref{lem4.13}:} Assume that $t-t_{0}%
\geq10r^{2}$ so that $\omega^{\ast}=\omega$. By Theorem \ref{lem4.10} we have
\begin{equation}
\omega(x,t,\Delta(x_{0},t_{0},r))\leq c|B_{d}(x,2r)|G(x,t,A_{2r}^{-}%
(x_{0},t_{0})), \label{lse1}%
\end{equation}
and using Lemma \ref{lem4.10b}, we get
\begin{equation}
|B_{d}(x,2r)|G(x,t,A_{2r}^{-}(x_{0},t_{0}))\leq c|B_{d}(x,2r)|G(x,t,A_{2r}%
^{+}(x_{0},t_{0})).
\end{equation}
Using (\ref{axel}) and the interior Harnack inequality (Lemma
\ref{lemmaacchar}) for the adjoint equation we get
\begin{equation}
|B_{d}(x,2r)|G(x,t,A_{2r}^{+}(x_{0},t_{0}))\leq c|B_{d}(x,r)|G(x,t,A_{r}%
^{+}(x_{0},t_{0})),
\end{equation}
where the constant is independent of $r$. Now, by Theorem \ref{lem4.10}, we
have
\begin{equation}
|B_{d}(x,r)|G(x,t,A_{r}^{+}(x_{0},t_{0}))\leq c\omega(x,t,\Delta(x_{0}%
,t_{0},r/2)). \label{lse4}%
\end{equation}
Combining (\ref{lse1})-(\ref{lse4}), we get the desired result when
$\omega^{\ast}=\omega.$ An analogous argument proves the theorem for
$\omega^{\ast}=\hat{\omega}$.\hfill\hfill$\Box$

\section{H\"{o}lder continuity of quotients\label{SecH}}

First we will prove a lemma, which is a generalization of Lemma 2 in
\cite{CG}. It is needed to get around the fact that we cannot approximate
$\Omega_{T}\cap C_{r}(x_{0},t_{0})$ with an NTA domain. In particular, the
approximation in Lemma \ref{regular} does not satisfy the inner corkscrew
condition, only the outer one.

For $(x_{0},t_{0}) \in S_{T}$, we let $U_{r}$ be a set which is regular for
the Dirichlet problem and which satisfies $C_{8r}(x_{0},t_{0}) \subset U_{r}
\subset C_{9r}(x_{0},t_{0})$. We let $\omega_{r}$ be the $X$-parabolic measure
with respect to $U_{r} \cap\Omega_{T}$, and $G_{r}$ be the Green function with
respect to $U_{r} \cap\Omega_{T}$.

Now set
\[
F=\partial_{p}\Omega_{T}\cap\left[  (B_{d}(x_{0},6r)\times\lbrack
t_{0}-36r^{2},t_{0}+4r^{2}])\setminus C_{2r}(x_{0},t_{0})\right]  .
\]
We can choose $(x_{i},t_{i})\in S_{T}$ such that $\{C_{\frac{r}{l}}%
(x_{i},t_{i})\}_{i=1}^{N}$ is a covering for $\overline{F}$, and
$\{C_{\frac{r}{100l}}(x_{i},t_{i})\}_{i=1}^{N}$ are disjoint. It is clear that
$N$ is independent of $r$ due to the definition of $F$. Now, we let
\[
H(x,t)=\sum_{i=1}^{N}\omega_{r}(x,t,\Delta_{2r/l}(x_{i},t_{i}))+|B_{d}%
(x,r)|G_{r}(x,t,A_{4r}^{-}(x_{0},t_{0})),
\]
and we note that $H(x,t)$ is a solution to the subelliptic heat equation. We
continue with a lemma and a theorem which constitutes the main tools in the
proof of H\"{o}lder continuity of quotients, Theorem \ref{holder}.

\begin{lemma}
\label{zombie} If $l$ is large enough, there exists $c=c(M)>0$ such that
$H(x,t)\geq c$ for $(x,t)\in\Omega_{T}\cap\partial_{p}C_{4r}(x_{0},t_{0})$.
\end{lemma}

\noindent\textbf{Proof:} By Lemma \ref{doftticka}, we have $\omega
_{r}(x,t,\Delta_{2r/l}(x_{i},t_{i}))\geq c$ for $(x,t)\in\Omega_{T}\cap
C_{r/l}(x_{i},t_{i})$. This proves the lemma in the case
\[
(x,t)\in\Omega_{T}\cap\partial_{p}C_{4r}(x_{0},t_{0})\cap\bigcup_{i=1}%
^{N}C_{r/l}(x_{i},t_{i}).
\]
Now take $(x,t)\in\lbrack\Omega_{T}\cap\partial_{p}C_{4r}(x_{0},t_{0}%
)]\setminus\bigcup_{i=1}^{N}C_{r/l}(x_{i},t_{i})$. It is easy to see that in
that case, $d(x,t,\partial_{p}\Omega_{T})>r/2l$. By first the Harnack
inequality and then Lemma \ref{lem4.10b}, we have
\begin{multline*}
|B_{d}(x,r)|G_{r}(x,t,A_{4r}^{-}(x_{0},t_{0}))\geq c|B_{d}(x,r)|G_{r}%
(A_{2r}^{-}(x_{0},t_{0}),A_{4r}^{-}(x_{0},t_{0}))\geq\\
\geq c|B_{d}(x,r)|G_{r}(A_{2r}^{+}(x_{0},t_{0}),A_{4r}^{-}(x_{0},t_{0}%
))=(\ast).
\end{multline*}
By the Harnack inequality, then Theorem \ref{lem4.10}, the Harnack inequality
again, and finally Lemma \ref{doftticka}, we get
\begin{multline*}
(\ast)\geq c|B_{d}(x,r)|G_{r}(A_{2r}^{+}(x_{0},t_{0}),A_{r}^{-}(x_{0}%
,t_{0}))\geq c\omega_{r}(A_{2r}^{+}(x_{0},t_{0}),\Delta(x_{0},t_{0}%
,r/2))\geq\\
\geq c\omega_{r}(A_{r/4}(x_{0},t_{0}),\Delta(x_{0},t_{0},r/2))\geq c.
\end{multline*}
\qquad We then have $H(x,t)\geq c>0$. Note that the boundary of $U_{r}$ does
not satisfy the inner corkscrew condition, however, none of the tools used in
this proof are such that they require this condition.\hfill$\Box$ \newline

\begin{remark}
\label{GULLVIVA} Let $\hat{\omega}_{r}$ be the adjoint $X$-parabolic measure
with respect to $U_{r}\cap\Omega_{T}$ and define
\[
F=\partial_{p}\Omega_{T}\cap\left[  (B_{d}(x_{0},6r)\times\lbrack t_{0}%
-4r^{2},t_{0}+36r^{2}])\setminus C_{2r}(x_{0},t_{0})\right]  .
\]
Then again we can choose $(x_{i},t_{i})\in S_{T}$ such that $\{C_{\frac{r}{l}%
}(x_{i},t_{i})\}_{i=1}^{N}$ is a covering for $\overline{F}$, and
$\{C_{\frac{r}{100l}}(x_{i},t_{i})\}_{i=1}^{N}$ are disjoint. Then, if we
define%
\[
\hat{H}(x,t)=\sum_{i=1}^{N}\hat{\omega}_{r}(x,t,\Delta_{2r/l}(x_{i}%
,t_{i}))+|B_{d}(x,r)|G_{r}(A_{4r}^{+},x,t),
\]
and if $l$ is large enough, there exists $c=c(M)>0$ such that $\hat
{H}(x,t)\geq c$ for $(x,t)\in\Omega_{T}\cap\partial_{p}C_{4r}(x_{0},t_{0}).$
This can be proved arguing as in Lemma \ref{zombie}, using properties of the
adjoint operator rather than properties of the operator itself. We omit the details.
\end{remark}

\begin{theorem}
\label{lem4.11} Assume (\ref{ass}). Let $u$, $v$ be non-negative solutions to
either the subelliptic heat equation or the adjoint subelliptic heat equation
in $\Omega_{T}\cap C_{20r}(x_{0},t_{0})$. Let $(x_{0},t_{0})\in S_{T}$ and let
$r<\min\{r_{0}/20,\sqrt{(T-t_{0})/20},\sqrt{t_{0}/20}\}$. Assume that $u$ and
$v$ vanish continuously on $\Delta(x_{0},t_{0},10r)$. Then there exists a
constant $c=c(H,M,\text{diam}(\Omega_{T}),T,r_{0})$, $c\geq1$, such that
\[
\frac{u(x,t)}{v(x,t)}\leq c\frac{u(A_{5r}^{+}(x_{0},t_{0}))}{v(A_{5r}%
^{-}(x_{0},t_{0}))}%
\]
if $u,v$ solves the subelliptic heat equation, and
\[
\frac{u(x,t)}{v(x,t)}\leq c\frac{u(A_{5r}^{-}(x_{0},t_{0}))}{v(A_{5r}%
^{+}(x_{0},t_{0}))}%
\]
if $u,v$ solves the adjoint subelliptic heat equation, whenever $(x,t)\in
\Omega_{T}\cap C_{r}(x_{0},t_{0})$.
\end{theorem}

\noindent\textbf{Proof:} Assume that $u$ and $v$ solves the subelliptic heat
equation. According to Lemma \ref{lem4.6} and Lemma \ref{zombie} we have
\begin{equation}
u(x,t)\leq Cu(A_{5r}^{+}(x_{0},t_{0}))H(x,t) \label{masher}%
\end{equation}
for $(x,t)\in\Omega_{T}\cap\partial_{p}C_{4r}(x_{0},t_{0})$. Let
\[
H^{\ast}(x,t)=\sum_{i=1}^{N}\omega_{r}(x,t,\Delta_{2r/l}(x_{i},t_{i}%
))+|B_{d}(x_{0},r_{0})|G_{r}(x,t,A_{4r}^{-}(x_{0},t_{0})).
\]
We introduce this function because, unlike $H(x,t)$, it will be a solution to
the subelliptic heat equation. By (\ref{axel}) we have $C^{-1}H(x,t)\leq
H^{\ast}(x,t)\leq CH(x,t)$ for $(x,t)\in\Omega_{T}\cap\partial_{p}C_{4r}%
(x_{0},t_{0})$. By the maximum principle we can then conclude that
(\ref{masher}) holds in $\Omega_{T}\cap C_{4r}(x_{0},t_{0})$.

We want to show that
\begin{equation}
\label{key}v(x,t) \geq Cu(A_{5r}^{-}(x_{0},t_{0})) H(x,t)
\end{equation}
for $(x,t) \in\Omega_{T} \cap C_{r}(x_{0},t_{0})$ (in that case we are done).

We can find $\alpha$ independent of $r$ such that $C_{\alpha}(A_{4r}^{-}%
(x_{0},t_{0}))\subset U_{r}\cap\Omega_{T}$. As in the proof of Theorem
\ref{lem4.10}, we let
\[
A=\{(x,t)\in U_{r}\cap\Omega_{T}:t=t_{0}-2(4r)^{2}\}\setminus C_{\alpha
r/2}(A_{4r}^{-}(x_{0},t_{0}))
\]%
\[
B=\{(x,t)\in\partial C_{\alpha r/2}(A_{4r}^{-}(x_{0},t_{0})):t>t_{0}%
-2(4r)^{2}\}.
\]
If $(x,t)\in A$, or $(x,t)\in\partial(U_{r}\cap\Omega_{T})$ with
$t>t_{0}-2(4r)^{2}$, then $|B_{d}(x,r)|G_{r}(x,t,A_{4r}^{-}(x_{0},t_{0}))=0$.
If $(x,t)\in B$, then Lemma \ref{lem4.3} and (\ref{axel}) imply that
$|B_{d}(x,r)|G_{r}(x,t,A_{4r}^{-}(x_{0},t_{0}))\leq C$. By the Harnack
principle, we have $v(A_{5r}^{-}(x_{0},t_{0}))\leq Cv(x,t)$ in $A$, so the
maximum principle gives
\begin{equation}
v(x,t)\geq C|B_{d}(x,r)|G_{r}(x,t,A_{4r}^{-}(x_{0},t_{0}))v(A_{5r}^{-}%
(x_{0},t_{0})) \label{imposition}%
\end{equation}
for $(x,t)\in\{(x,t)\in U_{r}\cap\Omega_{T}:t>t_{0}-2(4r)^{2}\}\setminus
C_{\alpha r/2}(A_{4r}^{-}(x_{0},t_{0})$.

We would now like to show that
\begin{equation}
H(x,t)\leq|B_{d}(x,r)|G_{r}(x,t,A_{4r}^{-}(x_{0},t_{0})) \label{attention}%
\end{equation}
in $\Omega_{T}\cap C_{r}(x_{0},t_{0})$. Note that $\Omega_{T}\cap C_{r}%
(x_{0},t_{0})\subset\Omega_{T}\setminus\bigcup_{i=1}^{N}C_{16r/l}(x_{i}%
,t_{i})$ if $l$ is small enough, so that we can apply Remark \ref{snippet} to
show that
\begin{equation}
\omega(x,t,\Delta_{2r/l}(x_{i},t_{i}))\leq C|B_{d}(x,{4r/l})|G_{r}%
(x,t,A_{4r/l}^{-}(x_{0},t_{0}))\leq C|B_{d}(x,r)|G_{r}(x,t,A_{4r}^{-}%
(x_{0},t_{0})) \label{revision}%
\end{equation}
when $(x_{i},t_{i})\in B_{d}(x_{0},6r)\setminus B_{d}(x_{0},2r)$. When
$(x_{i},t_{i})\in B_{d}(x_{0},2r)\times\lbrack t_{0}-36r^{2},t_{0}-4r^{2}]$,
we can get (\ref{revision}) by means of Theorem \ref{lem4.10} instead, which
completes the proof of (\ref{attention}). By (\ref{attention}) and
(\ref{imposition}) combined, we get (\ref{key}).

If $u$ and $v$ solve the adjoint subelliptic heat equation, we can argue in a
similar manner, using Remark \ref{GULLVIVA} instead of Lemma \ref{zombie} and
with $H^{\ast}(x,t)$ replaced by
\[
\hat{H}^{\ast}(x,t)=\sum_{i=1}^{N}\hat{\omega}_{r}(x,t,\Delta_{2r/l}%
(x_{i},t_{i}))+|B_{d}(x_{0},r_{0})|G_{r}(A_{4r}^{+}(x_{0},t_{0}),x,t),
\]
which is a solution to the adjoint subelliptic heat equation. We omit the
details. \hfill$\Box$ \newline

We now have the necessary tools to prove Theorem \ref{holder}. In fact, the
proof follows directly using the proof of Lemma 4.2 and Theorem 1.2 in
\cite{FGGMN} when $u,v$ solves the subelliptic heat equation. If $u,v$ solves
the adjoint subelliptic equation, the proof is similar, using the adjoint
estimates. We omit the details, but remark that when we refer to Lemma 3.3 in
\cite{FGGMN}, this corresponds to Lemma \ref{lem4.5} in the present paper.

\section{Fatou Theorem\label{SecFatou}}

\label{beshem}

The goal of this section is to prove Theorem \ref{fatou}. Our approach rely on
the use of certain kernel functions for $H$ in $\Omega$, and for this purpose,
let $S$ denote the lateral boundary of $\Omega$.

\begin{definition}
\label{flapjack} A function $K:\Omega\times\partial_{p}\Omega\rightarrow
\mathbb{R}^{+}\cup\{+\infty\}$ is a kernel function at $(y,s)\in\partial
_{p}\Omega$, normalized in $(X_{0},T_{1})\in\Omega$, if \newline(i)
$K(X_{0},T_{1},y,s)=1$,\newline(ii) $K(\cdot,\cdot,y,s)\geq0$ is a weak
solution of $Hu=0$ in $\Omega$.\newline(iii) $K(\cdot,\cdot,y,s)$ is
continuous in $\overline{\Omega}\setminus\{(y,s)\}$ and
\[
\lim_{(x,t,)\rightarrow(y_{0},s_{0})}K(x,t,y,s)=0
\]
if $(y_{0},s_{0})\in\partial_{p}\Omega_{T}\setminus\{(y,s)\}$. If $s>T$, we
will let $K(x,t,y,s)=0$.
\end{definition}

\begin{lemma}
\label{hjortron} Assume (\ref{ass}). Then there exists a unique kernel
function (for $H$ and $\Omega$) at $(y,s) \in\partial_{p} \Omega$ normalized
for $(X_{0},T_{1})$ if $0 \leq s \leq T_{1}- \delta^{2}$.
\end{lemma}

\noindent\textbf{Proof:} Let $(y,s)\in\partial_{p}\Omega$ with $0\leq s\leq
T_{1}-\delta^{2}$ and define the following
\begin{align*}
A_{n}^{+}  &  =A_{2^{-n}}^{+}(y,s),\\
\Phi_{n}  &  =\{(x,t):d(x,y)<2^{-n},-4^{-n}<s-t<2\cdot4^{-n}\},\\
\Omega_{T_{1}}^{n}  &  =(\Omega_{T_{1}}\cap\{t>s\})\setminus\overline{\Phi
_{n}},\\
\alpha_{n}  &  =\partial_{p}\Phi_{n}\cap\overline{\Omega_{T_{1}}^{n}},\\
\beta_{n}  &  =\{(x,t)\in\alpha_{n}:t>s+4^{-n},d(x,\partial\Omega
)>2^{-n}/2M\}.
\end{align*}

Note that $\beta_{n}$ is not empty, since $A_{n}^{+}(y,s) \in\beta_{n}$. We
let $\omega_{n}$ be the $X$-caloric measure with regard to the domain
$\Omega_{T_{1}}^{n}$ (note that this is possible since $\Omega_{T_{1}}^{n}$ is
regular for the Dirichlet problem; see Lemma \ref{regular}).

As in the proof of Theorem 2.7 in \cite{FGS}, we prove that
\begin{equation}
\label{bronte}K(x,t,y,s) = \lim_{n \to\infty} \frac{\omega_{n}(x,t,\beta_{n}%
)}{\omega_{n}(X_{0},T_{1},\beta_{n})}%
\end{equation}
is a kernel function. To prove uniqueness, we first claim that if $v$ is
another kernel function at $(y,s)$, normalized at $(X_{0},T_{1})$, then there
exists a constant $C=C(H,M,r_{0})$ such that
\begin{equation}
\label{paus}v(x,t) \geq CK(x,t,y,s)
\end{equation}
for all $(x,t) \in\Omega$ (note that it is trivially true for $t<s$). If this
claim is true, we proceed as in the proof of Theorem 1.7 in Kemper \cite{K}.

To prove (\ref{paus}), we first observe that
\[
v(x,t)=\int_{\partial_{p}\Omega_{T_{1}}^{n}}v(\xi,\tau)\omega_{n}(x,t,\xi
,\tau)\geq\inf_{\beta_{n}}v(x,t)\omega_{n}(x,t,\beta_{n}),
\]
and that due to the Harnack inequality, we have
\[
\inf_{\beta_{n}}v(x,t)\geq Cv(A_{n+1}^{+}).
\]
Combining these, we get
\begin{equation}
v(x,t)\geq Cv(A_{n+1}^{+})\omega_{n}(x,t,\beta_{n}) \label{crash}%
\end{equation}
for $(x,t)\in\Omega_{T_{1}}^{n}$. By Lemma \ref{smultron} we also have
$v(x,t)\leq Cv(A_{n+1}^{+})$ for $(x,t)\in\Omega\setminus C_{2^{-n}}(y,s)$.
The maximum principle gives
\begin{equation}
v(x,t)\leq Cv(A_{n+1}^{+})\omega_{n}(x,t,\alpha_{n}), \label{freeze}%
\end{equation}
for $(x,t)\in\Omega_{T_{1}}^{n}$, since $\alpha_{n}\subset\Omega\setminus
C_{2^{-n}}(y,s)$. Combining (\ref{crash}) and (\ref{freeze}), and setting
$(x,t)=(X_{0},T_{1})$ in the latter, we get
\[
v(x,t)\geq C\frac{\omega_{n}(x,t,\beta_{n})}{\omega_{n}(X_{0},T_{1},\alpha
_{n})}%
\]
for $(x,t)\in\Omega_{T_{1}}^{n}$. To prove (\ref{paus}), it is enough to prove
that
\begin{equation}
\omega_{n}(X_{0},T_{1},\alpha_{n})\leq C\omega_{n}(X_{0},T_{1},\beta_{n}).
\label{southpaw}%
\end{equation}
Since $\Omega$ is an NTA domain, we can find a point $(x_{1},t_{1})\in
\beta_{n}$ such that $\Delta_{2^{-n}/2M}(x_{1},t_{1})\subset\beta_{n}$ and
such that $t_{1}=t_{0}+2\cdot4^{-n}$, that is, such that $(x_{1},t_{1})$ is on
the \textquotedblleft top\textquotedblright\ of $\beta_{n}$. By the Harnack
inequality and the proof of Lemma \ref{doftticka}, we have
\[
\omega_{n}(X_{0},T_{1},\beta_{n})\geq\omega_{n}(X_{0},T_{1},\Delta_{2^{-n}%
/2M}(x_{1},t_{1}))\geq C\omega_{n}(x_{1},t_{1}+2^{-2n}/8M^{2},\Delta
_{2^{-n}/2M}(x_{1},t_{1}))>C.
\]
Since $\omega_{n}(X_{0},T_{1},\alpha_{n})\leq1$, we have proved
(\ref{southpaw}), so we are done. \hfill$\Box$ \newline

Now that we know that the kernel is unique, we can also define it as
\begin{equation}
K(x,t,y,s)=\lim_{\epsilon\rightarrow\infty}\frac{\omega(x,t,\Delta_{\epsilon
})}{\omega(X_{0},T_{1},\Delta_{\epsilon})}, \label{shenandoah}%
\end{equation}
since this expression also satisfies Definition \ref{flapjack} in the same way
that (\ref{bronte}) does. From now on, when working in $\Omega_{T}$ we will
always assume that the kernel is normalized in $T_{1}=T+1$ so as to avoid the
limitation $t<T_{1}-\delta^{2}$ in the previous theorem.

\begin{lemma}
\label{lastbil} Assume (\ref{ass}). Let $(x_{0},t_{0})\in S_{T}$. If
$r<\min\{r_{0}/2,\sqrt{T-t_{0}}/4,\sqrt{t_{0}}/4\}$, then
\[
\lim_{(x,t)\rightarrow(x_{0},t_{0})}\sup\{K(x,t,y,s):(y,s)\in\partial
_{p}\Omega_{T}\setminus\Delta(x_{0},t_{0},r)\}=0.
\]

\end{lemma}

\noindent\textbf{Proof:} Let $(x_{0},t_{0})\in S_{T}$. Using Lemma
\ref{regular}, we let $D_{r}$ be a set which is regular for the Dirichlet
problem and such that $B_{d}(x_{0},r/8)\subset D_{r}\subset B_{d}(x_{0},r/4)$.
Let $\tilde{D}_{r}=D_{r}\times(t_{0}-r^{2}/16,t_{0}-r^{2}/16)$ and let
$h_{r}(x,t)$ be the H-parabolic measure of $\partial_{p}\tilde{D}_{r}%
\cap\Omega_{T}$ relative to $\tilde{D}_{r}\cap\Omega_{T}$. By Lemma
\ref{lem4.6} and the Harnack inequality, we have
\[
K(x,t,y,s)\leq CK(A_{r/4}^{+}(x_{0},t_{0}),y,s)\leq CK(X_{0},T_{1},y,s)=C
\]
when $(x,t)\in C_{r/4}(x_{0},t_{0})\cap\Omega_{T}$. Note that $C$ in fact
depend on $r$. Using the maximum principle and then taking the supremum over
$(y,s)$, we get
\[
\sup\{K(x,t,y,s):(y,s)\in\partial_{p}\Omega_{T}\setminus\Delta(x_{0}%
,t_{0},r)\}\leq Ch_{r}(x,t)
\]
for $(x,t)\in\tilde{D}_{r}\cap\Omega_{T}$. Letting $(x,t)\rightarrow
(x_{0},t_{0})$, we are finished, since $h_{r}(x_{0},t_{0})=0$. \hfill$\Box$
\newline

\begin{remark}
We can also prove Lemma \ref{lastbil} when $t=0$, using \textit{(iii) }in the
definition of a kernel.
\end{remark}

We can prove the following lemma as in the proof of Theorem 1.10 in \cite{K},
using Lemma \ref{tag} in the present paper when they refer to Lemma 1.4 in
\cite{K}.

\begin{lemma}
\label{mint} Assume (\ref{ass}). If $u$ is a non-negative solution of $Hu=0$
in $\Omega_{T}$, there exists a Borel measure $\nu,$ depending on $u$, on
$\partial_{p}\Omega_{T}$, such that
\[
u(x,t)=\int_{\partial_{p}\Omega_{T}}K(x,t,y,s)d\nu(y,s),
\]
where $K(x,t,y,s)$ is the kernel function for $H$ and $\Omega_{T}$, normalized
at $(X_{0},T_{1})$.
\end{lemma}

\begin{lemma}
\label{kantarell} Assume (\ref{ass}). Let $r<\min\{r_{0}/2,\sqrt{(T-t_{0}%
)/4},\sqrt{t_{0}/4}\}$, $(x_{0},t_{0})\in S_{T}$. Let $u$ and $v$ be two
non-negative solutions of $Hu=0$ in $\Omega_{T_{1}}$, and assume that $u$ and
$v$ vanish continuously on $\partial_{p}\Omega_{T_{1}}\setminus\Delta
_{r/2}(x_{0},t_{0})$. Then there exists a constant $c=c(H,M,r_{0}%
,X_{0},diam(\Omega),T)$ such that
\[
u(X_{0},T_{1})v(A_{r}^{+}(x_{0},t_{0}))\leq Cv(X_{0},T_{1})u(A_{r}^{+}%
(x_{0},t_{0})).
\]

\end{lemma}

\noindent\textbf{Proof:} By Lemma \ref{tag} we have $u(X_{0},T_{1})\leq
Cu(A_{r}^{+}(x_{0},t_{0}))\omega(X_{0},T_{1},\Delta_{2r}(x_{0},t_{0}))$, for
$(X_{0}T_{1})\in\Omega_{T}\backslash C_{r}(x_{0},t_{0})$. Let $\Omega
_{r}=\{(x,t)\in\Omega_{T_{1}}:t>s\}\setminus\Phi_{r}$, where $\Phi_{r}$ is the
same as $\Phi_{n}$ in the proof of Lemma \ref{hjortron}, but with $2^{-n}$
replaced with $r$. We then get
\[
v(X_{0},T_{1})\geq Cv(A_{r}^{+}(x_{0},t_{0}))\omega_{n}(X_{0},T_{1}%
,\beta_{r/2}),
\]
in the same way as (\ref{crash}). Combining what we have so far, we get
\[
u(X_{0},T_{1})v(A_{r}^{+}(x_{0},t_{0}))\leq C\frac{v(X_{0},T_{1})u(A_{r}%
^{+}(x_{0},t_{0}))\omega(X_{0},T_{1},\Delta_{2r}(x_{0},t_{0}))}{\omega
_{n}(X_{0},T_{1},\beta_{r/2})}.
\]
Finally, we argue as in the proof of (\ref{southpaw}) to obtain
\[
\frac{\omega(X_{0},T_{1},\Delta_{2r}(x_{0},t_{0}))}{\omega_{n}(X_{0}%
,T_{1},\beta_{r/2})}\leq C,
\]
which completes the proof. \hfill$\Box$ \newline

\begin{lemma}
\label{nektarin} Assume (\ref{ass}). Let $(x_{0},t_{0}) \in\partial_{p}
\Omega_{T}$. If $r < \min\{r_{0}/2,\sqrt{(T-t_{0})/4},\sqrt{t_{0}/4}\}$, there
exists a constant $c = c(H,M,r_{0},X_{0},diam(\Omega),T)$ such that
\[
\sup_{(y,s) \in\Delta(x_{0},t_{0},r)} K(A^{+}_{2r}(x_{0},t_{0}),y,s) \leq
\frac{C}{\omega_{X}(X_{0},T_{1}, \Delta(x_{0},t_{0},r))}.
\]

\end{lemma}

\noindent\textbf{Proof:} The proof is essentially the same as for Theorem 2.11
in \cite{FGS}. It uses Lemma \ref{kantarell} and (\ref{shenandoah}).
\hfill$\Box$ \newline

For $(y,s)\in\partial_{p}\Omega_{T}$, we let $\Delta_{j}(y,s)=\Delta
(y,s,2^{j}r)$ and $R_{j}(y,s)=\Delta_{j}(y,s)\setminus\Delta_{j-1}(y,s)$.

\begin{lemma}
\label{kladdkaka} Assume (\ref{ass}). Let $(x_{0},t_{0})\in\partial_{p}%
\Omega_{T}$ and let $r<\min\{r_{0}/2,\sqrt{(T-t_{0})/4},\sqrt{t_{0}/4}\}$.
Then there exists a sequence $\{C_{j}\}$ of positive numbers, which depend on
$H$, $M$, $r_{0}$, $X_{0}$, $diam(\Omega)$, and $T$ (but not on $r$ or
$(x_{0},t_{0})$), such that
\begin{equation}
\sup_{(y,s)\in R_{j}(x_{0},t_{0})}K(A_{2r}^{+}(x_{0},t_{0}),y,s)\leq
\frac{C_{j}}{\omega(X_{0},T_{1},\Delta_{j}(x_{0},t_{0}))}. \label{isolation}%
\end{equation}
Moreover, $\sum_{j}C_{j}<\infty$.
\end{lemma}

\noindent\textbf{Proof: }As a first step we prove that (\ref{isolation}) holds
for a constant $C$ independent of our choice of $j$. Now, for $j=1,...,8$, we
have;%
\begin{align*}
\sup_{(y,s)\in R_{j}(x_{0},t_{0})}K(A_{2r}^{+}(x_{0},t_{0}),y,s)  &  \leq
\sup_{(y,s)\in\Delta_{j}(x_{0},t_{0})}K(A_{2r}^{+}(x_{0},t_{0}),y,s)\\
&  \leq c\sup_{(y,s)\in\Delta_{j}(x_{0},t_{0})}K(A_{2^{j+1}r}^{+}(x_{0}%
,t_{0}),y,s)\\
&  \leq\frac{C}{\omega(X_{0},T_{1},\Delta_{j}(x_{0},t_{0}))}.
\end{align*}
Above we used the definition of $R_{j}$ and $\Delta_{j}$ in the first
inequality, the Harnack principle in the second inequality and Lemma
\ref{nektarin} in the last inequality. For $j>8$, we know by Lemma
\ref{nektarin} that,
\begin{equation}
\sup_{(y,s)\in R_{j}(x_{0},t_{0})}K(A_{2^{j+1}r}^{+}(x_{0},t_{0}%
),y,s)\leq\frac{C}{\omega(X_{0},T_{1},\Delta_{j}(x_{0},t_{0}))}.
\label{gatwick}%
\end{equation}
Now, we observe that for $(y,s)\in R_{j}(x_{0},t_{0})$, $K(\cdot,\cdot,y,s)$
is a non-negative solution to the subelliptic heat equation, vanishing on (at
least)%
\[
\partial_{p}\Omega_{T}\backslash\Delta_{2^{j-5}r}(y,s).
\]
Therefore, we use Lemma \ref{smultron} to get%
\begin{equation}
K(x,t,y,s)\leq CK(A_{2^{j-4}r}^{+}(y,s),y,s), \label{ettan}%
\end{equation}
for all $(x,t)\in\Omega_{T}\backslash C_{2^{j-4}r}(y,s)$. For $(y,s)\in
R_{j}(x_{0},t_{0}),$ let
\[
(y_{0,r},s_{0,r})=A_{2^{j+1}r}^{+}(x_{0},t_{0})\text{ \ \ and \ \ }%
(y_{r},s_{r})=A_{2^{j-4}r}^{+}(y,s).
\]
Then $|y_{0,r}-y_{r}|$ is bounded by $C2^{j}r$, while $s_{0,r}-s_{r}%
>2^{2j+1}r^{2}$. Hence, using the Harnack principle we obtain%
\begin{equation}
K(A_{2^{j-4}r}^{+}(y,s),y,s)\leq CK(A_{2^{j+1}r}^{+}(x_{0},t_{0}),y,s).
\label{tvaan}%
\end{equation}
Further, using (\ref{ettan}) and (\ref{tvaan}), and for all $(x,t)\in
\Omega_{T}\backslash C_{2^{j-4}r}(y,s),$ we have that%
\begin{equation}
K(x,t,y,s)\leq CK(A_{2^{j+1}r}^{+}(x_{0},t_{0}),y,s). \label{maskar}%
\end{equation}
Let $(x,t)=A_{2r}^{+}(x_{0},t_{0})$ in (\ref{maskar}), and using
(\ref{gatwick}), we get%
\[
\sup_{(y,s)\in R_{j}(x_{0},t_{0})}K(A_{2r}^{+}(x_{0},t_{0}),y,s)\leq\frac
{C}{\omega(X_{0},T_{1},\Delta_{j}(x_{0},t_{0})}.
\]
This is indeed justified, since $A_{2r}^{+}(x_{0},t_{0})\in\Omega
_{T}\backslash C_{2^{j-4}r}(y,s)$ for $j>8$.

By Lemma \ref{regular}, we can let $U_j$ be a set which is regular for the Dirichlet problem and which satisfies $B_d(x_0,2^{j-2}r) \subseteq U_j \subseteq B_d(x_0,2^{j-1}r)$. Let $\Sigma_{j}= \Omega_T\cap (U_j \times [t_0-4^{j-1}r, t_0+4^{j-1}r])$, and let $h_{j}$ be the $X$-caloric measure of $\Omega
_{T}\cap\partial_{p}\Sigma_{j}$ with respect to $\Sigma_{j},$ that is,
\[
h_{j}(x,t)=\omega(x,t,\Omega_{T}\cap\partial_{p}\Sigma_{j},\Sigma_{j})\text{.}
\]
By definition, we have the following;%
\[
h_{j}(A_{2r}^{+}(x_{0},t_{0}))=\omega(A_{2r}^{+}(x_{0},t_{0}),\Omega_{T}
\cap\partial_{p}\Sigma_{j},\Sigma_{j}),
\]
and, in particular, $h_{j}$ is the solution to%
\[
\left\{
\begin{array}
[c]{ll}%
Hh_{j}=0 & \text{in }\Sigma_{j},\\
h_{j}=\chi_{\lbrack\Omega_{T}\cap\partial_{p}\Sigma_{j}]} & \text{on }%
\partial_{p}\Sigma_{j}.
\end{array}
\right.
\]
By the maximum principle, we then have%
\[
\sup_{(y,s)\in R_{j}(x_{0},t_{0})}K(A_{2r}^{+}(x_{0},t_{0}),y,s)\leq
\frac{Ch_{j}(A_{2r}^{+}(x_{0},t_{0}))}{\omega(X_{0},T_{1},\Delta_{j}%
(x_{0},t_{0}))}.
\]
What is left to prove is that $\sum_{j}h_{j}(A_{2r}^{+}(x_{0},t_{0}))<\infty$,
and we will use Lemma \ref{lem4.5} to do that. First, we note that
\[
d_{p}(A_{2r}^{+}(x_{0},t_{0}),x_{0},t_{0}) \leq ((2r)^2+2(2r)^2))^{1/2} \leq \sqrt{12} r \leq 4r
\]
The largest cylinder $C_{4R}$ we can consider when using Lemma \ref{lem4.5} has radius $4R= 2^{j-2}r$, or $R=2^{j-4}r$. Further, $A_{2r}^{+}(x_{0},t_{0})\in C_{R}(x_{0},t_{0})$ if $R^{2}%
>2(2r)^{2}=8r^{2}$, that is, $R>2\sqrt{2}r$. In particular, the
calculations that follow apply for all $j>8$. Now, using Lemma \ref{lem4.5}
on $\Sigma_{j}$, we get
\begin{align}
h_{j}(A_{2r}^{+}(x_{0},t_{0}))  &  \leq c\left(  \frac{4r}{2^{j-6}r}\right)  ^{\alpha}
\sup_{\Omega_{T}\cap C_{2^{j-5}r}(x_{0},t_{0})}
h_{j}\nonumber\\
&  \leq c (2^{-j})^{\alpha}.
\label{nfl}%
\end{align}
$\allowbreak\allowbreak\allowbreak$ In the first inequality we used Lemma
\ref{lem4.5}, in the second that $h_{j}(x,t)\in\lbrack0,1]$. We know that
$\alpha\in(0,1)$, and therefore $2^{\alpha}>1$ which means that $\sum
_{j>8}(2^{\alpha})^{-j}$ converges. In particular, using this observation with
(\ref{nfl}), and the fact that $h_{j}(A_{2r}^{+}(x_{0},t_{0}))$ is finite for
$j=1,...,8,$ we obtain
\[
\sum_{j=1}^{\infty}h_{j}(A_{2r}^{+}(x_{0},t_{0}))<\infty.
\]
This completes the proof. \hfill$\Box$ \newline

We define a non-tangential region at $P \in\partial_{p} \Omega_{T}$ as
\[
\Gamma_{\alpha}(P) = \{ (x,t) \in\Omega_{T}: d_{p}(x,t,P) \leq(1+ \alpha)
d_{p}(x,t,\partial_{p} \Omega_{T}) \},
\]
and the non-tangential maximal function as $N_{\alpha}(u)(P) = \sup_{(x,t)
\in\Gamma_{\alpha}(P)} |u(x,t)|$. We also define the Hardy-Littlewood maximal
function of the measure $\nu$ with respect to $\omega(X_{0},T_{1})$ as
\[
M_{\omega}(\nu)(x,t) = \sup_{r>0} \frac{\nu(\Delta(x,t,r))}{\omega(X_{0}%
,T_{1},\Delta(x,t,r))}.
\]

\begin{lemma}
\label{fika} Assume (\ref{ass}). If $\nu$ is a finite Borel measure on
$\partial_{p} \Omega_{T}$ such that
\[
u(x,t) =\int_{\partial_{p} \Omega_{T}} K(x,t,y,s) d \nu(y,s),
\]
then for every $(y,s) \in\partial_{p} \Omega_{T}$ there exists a constant $c =
c(H,M,r_{0},diam(\Omega_{T}),T,\alpha,\omega_{X}(X_{0},T_{1}))$ such that
\[
N_{\alpha}(u) (y,s) \leq c M_{\omega}(\nu) (y,s).
\]

\end{lemma}

\noindent\textbf{Proof:} The proof is essentially the same as Theorem 2.13 in
\cite{FGS}. It uses Lemma \ref{mint}, Lemma \ref{kladdkaka}, the Harnack
inequality and Lemma \ref{lem4.9}. \hfill$\Box$ \newline

\noindent\textbf{Proof of Theorem \ref{fatou}:} By Lemma \ref{mint} we can
write
\[
u(x,t)=\int_{\partial_{p}\Omega}K(x,t,y,s)d\nu(y,s),
\]
and we can make the decomposition $d\nu=fd\omega(X_{0},T_{1})+d\nu_{s}$, where
$d\nu_{s}\perp d\omega(X_{0},T_{1})$ and $f\in L^{1}(\partial_{p}%
\Omega,d\omega(X_{0},T_{1}))$. So we have
\[
u(x,t)=\int_{\partial_{p}\Omega_{T}}f(y,s)K(x,t,y,s)d\omega(X_{0}%
,T_{1},y,s)+\int_{\partial_{p}\Omega}K(x,t,y,s)d\nu_{s}(y,s)=u_{a}%
(x,t)+u_{s}(x,t).
\]
Let $F$ be the support of $\nu_{s}$. A standard strategy using Lemma
\ref{fika} (see for instance \cite[Theorem 4.4]{Caf}) shows that
$u_{a}(x,t)\rightarrow f(y,s)$ when $(x,t)\rightarrow(y,s)$ in $\Gamma
_{\alpha}(y,s)$ almost everywhere on $\partial_{p}\Omega$ with respect to
$\omega(X_{0},T_{1})$. Moreover, if $(y,s)\not \in F$, $u_{s}(x,t)\rightarrow
0$ when $(x,t)\rightarrow(y,s)$ in $\Gamma_{\alpha}(y,s)$ by Lemma
\ref{lastbil}. This concludes the proof.\hfill$\Box$

\subsection*{Acknowledgements}

The authors wish to thank Kaj Nystr\"{o}m for bringing the problem to their attention.

\end{document}